\documentclass[12pt]{iopart}

\pdfoutput=1 

\usepackage{graphicx}
\usepackage{color}
\usepackage{hyperref}
\usepackage{ifthen}
\usepackage{amsfonts} % needed by fancy expectation symbol

\begin{document}

% https://tex.stackexchange.com/questions/56765/getting-the-expectation-symbol-to-behave-like-sum-instead-of-sigma
% \DeclareMathOperator{\E}{\mathbb{E}}
\newcommand{\E}{\mathbb{E}}
\newcommand{\V}{\mathbb{V}}
\newcommand{\W}{\mathbb{W}}

\newcommand \completegraph[1][] {
	\ifthenelse{ \equal{#1}{} }
	{ \mathcal{K}_n }
	{ \mathcal{K}_{#1} }
}
\newcommand \lineartree[1][] {
	\ifthenelse{ \equal{#1}{} }
	{ \mathcal{L}_n }
	{ \mathcal{L}_{#1} }
}
\newcommand \startree[1][] {
	\ifthenelse{ \equal{#1}{} }
	{ \mathcal{S}_n }
	{ \mathcal{S}_{#1} }
}

\title{The sum of edge lengths in random linear arrangements}

\author{Ramon Ferrer-i-Cancho}
\address{Complexity \& Quantitative Linguistics Lab, LARCA Research Group \\
Departament de Ci\`encies de la Computaci\'o, \\
Universitat Polit\`ecnica de Catalunya, \\
Campus Nord, Edifici Omega, Jordi Girona Salgado 1-3. \\
08034 Barcelona, Catalonia (Spain)}
\ead{rferrericancho@cs.upc.edu}
\begin{abstract}
Spatial networks are networks where nodes are located in a space equipped with a metric. Typically, the space is two-dimensional and until recently and traditionally, the metric that was usually considered was the Euclidean distance. In spatial networks, the cost of a link depends on the edge length, i.e. the distance between the nodes that define the edge. Hypothesizing that there is pressure to reduce the length of the edges of a network requires a null model, e.g., a random layout of the vertices of the network. Here we investigate the properties of the distribution of the sum of edge lengths in random linear arrangement of vertices, that has many applications in different fields. A random linear arrangement consists of an ordering of the elements of the nodes of a network being all possible orderings equally likely. The distance between two vertices is one plus the number of intermediate vertices in the ordering. Compact formulae for the 1st and 2nd moments about zero as well as the variance of the sum of edge lengths are obtained for arbitrary graphs and trees. We also analyze the evolution of that variance in Erd\H{o}s-R\'enyi graphs and its scaling in uniformly random trees. Various developments and applications for future research are suggested.       
\end{abstract}

\noindent {\small {\it Keywords\/}: networks, random linear arrangement, scaling laws}

\pacs{89.75.Hc Networks and genealogical trees \\
89.75.Fb Structures and organization in complex systems \\
89.75.Da Systems obeying scaling laws}

% \submitto{\JSTAT}

\maketitle

\section{Introduction}

Spatial networks are networks for which the nodes are located in a space equipped with a metric \cite{Barthelemy2011a}.
For most practical applications, the space is two-dimensional and the metric is the usual Euclidean distance \cite{Barthelemy2011a}. 
Non-Euclidean spaces have been introduced in network research for the ease with which they reproduce the heterogeneous degree distributions that are found in real complex networks \cite{Krioukov2010a}. 

The length of an edge is defined as the metric distance between the nodes that form it (see for a discussion of the suitability of the term length \cite{Ferrer2017c}). A fundamental implication of space on networks is that links have a cost that depends on their length \cite{Barthelemy2011a}. A important example are brain networks, where regions that are spatially closer have a greater probability of being connected than remote
regions as longer axons are more costly in terms of material and energy \cite{Bullmore2009a}.
In general, to argue that there is pressure to reduce the length of edges a null hypothesis is necessary (e.g., \cite{Ferrer2004b}).

Many models of spatial networks, e.g., random geometric graphs, assume that both the structure of the network and the layout of vertices is random \cite{Barthelemy2018a}.  
Here we are interested in the particular problem of a network whose structure is given {\em a priori} and their nodes are arranged in a $1$-dimensional Euclidean space \cite{Ferrer2004b}. For simplicity, let us suppose that the vertices are arranged linearly, namely, forming a sequence. The $i$-th vertex of the sequence and the $j$-th vertex of the sequence are at distance $|i-j|$ and the length of an edge that joins them is then $|i-j|$. Put differently, the distance between two vertices is the number of intermediate vertices in the linear arrangement plus one. Here we aim to investigate the statistical properties of the distribution of edge lengths in random linear arrangements. 

Suppose a network of $n$ vertices and $m$ edges. 
The sum of edge lengths of a linear arrangement of the vertices of that network is defined as
\begin{equation}
D = \sum_{i=1}^m d_i,
\label{sum_of_edge_lengths_equation} 
\end{equation} 
where $d_i$ is the length of the $i$-th edge. Equivalently, it can be defined as
\begin{equation*}
D = \sum_{d=1}^{n-1} m(d) d,
\end{equation*} 
where $m(d)$ is the number of edges of length $d$.

\begin{figure}
\begin{center}
\includegraphics[scale = 0.6]{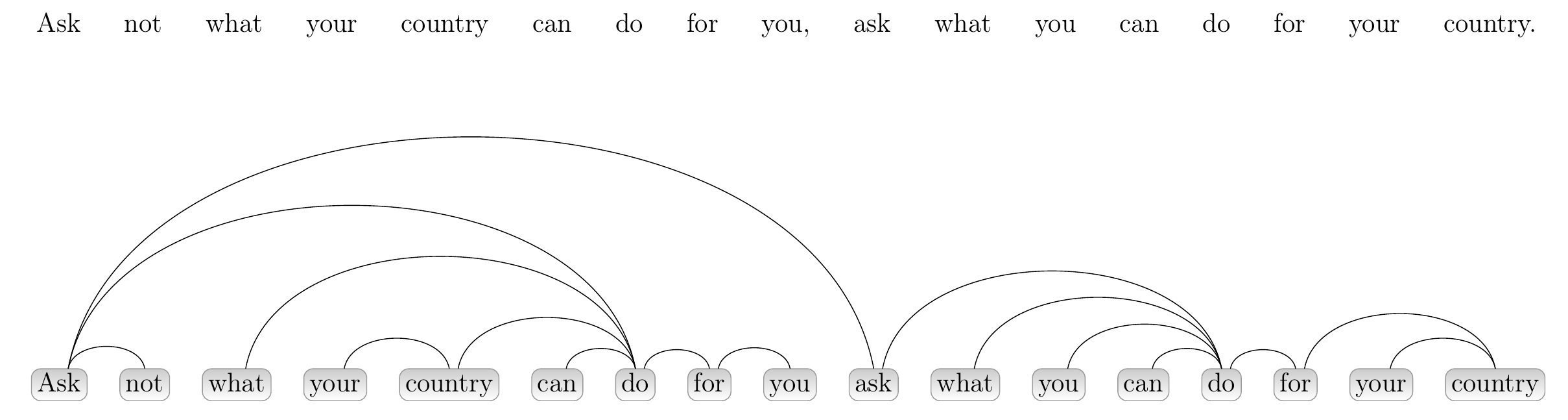}
\end{center}
\caption{\label{syntactic_dependency_tree_figure} The linear arrangement of the words of a sentence (top) and the corresponding network where vertices are words and edges indicate syntactic dependencies (bottom). Adapted from \url{https://cloud.google.com/natural-language/docs/morphology}. For simplicity, link directions are omitted. Punctuation marks are excluded as vertices following standards from research on dependency lengths \cite{Ferrer2016a}. }
\end{figure}
 
The statistical properties of linear arrangements of given networks are relevant in many contexts. In linguistics, the structure of a sentence can be defined as network where vertices are words and edges indicate syntactic dependencies between words (Fig. \ref{syntactic_dependency_tree_figure}). The linear order is defined by the sequential order of the words in the sentence. For the network in Figure \ref{syntactic_dependency_tree_figure}, $n = 17$, $m = 16$ and $D= 40$ (Table \ref{lengths_table}).
In these kind of networks, $D$ has been shown to be smaller than expected by chance, and pressure to reduce the distance between connected words is believed to result from two factors: decay of activation and interference \cite{Liu2017a}. In computer science, the minimum linear arrangement problem consists of finding $D_{min}$, the minimum value of $D$ over all possible $n!$ linear arrangements \cite{Diaz2002, Petit2011a}. $\E_{rlt}[D_{min}]$, the expectation of $D_{min}$ in the ensemble of uniformly random labeled trees (with vertex labels as vertex positions) grows logarithmically $n \geq 3$, i.e. \cite{Esteban2016a} 
\begin{equation}
\E_{rlt}[D_{min}/(n-1)] \approx a \log n + b,
\label{scaling_of_mean_mla_equation}
\end{equation}
where $a$ and $b$ are two constants and $D_{min}/(n-1)$ is the mean length of edges in a minimum linear arrangement of one of these trees (a tree has $n - 1$ edges).

\begin{table}
\caption{\label{lengths_table} $m(d)$, the number of edges of length $d$, for the network in Fig. \ref{syntactic_dependency_tree_figure}. $m(d)=0$ for $d > 9$. }
\begin{indented}
\item[]
\begin{tabular}{@{}rr}
\br
$d$ & $m(d)$ \\ % \hline
\mr
1 & 8 \\
2 & 3 \\
3 & 1 \\
4 & 2 \\
5 & 0 \\
6 & 1 \\
7 & 0 \\
8 & 0 \\
9 & 1 \\
\br
\end{tabular}
\end{indented}
\end{table}

Here we are interested in some properties of the distribution of $D$ in uniformly random linear arrangements ({\em rla}). In other words, we are interested in the statistical properties of $D$ over the ensemble of the $n!$ equally likely orderings of the vertices of a network. These random orderings provide a baseline or null model to study the properties of the actual linear arrangements of real networks \cite{Ferrer2004b,Gomez2016a}. 

Here we aim to calculate $\V_{rla}[D]$, the variance of $D$ in a uniformly random linear arrangement, that is defined as
\begin{equation*}
\V_{rla}[D] = \E_{rla}[D^2] - \E_{rla}[D]^2, 
\end{equation*}
where $\E_{rla}[D]$ and $\E_{rla}[D^2]$ are the 1st and the 2nd moment about zero of $D$, respectively. 
When the network has only one edge ($m=1$), $D$ matches $d$, the length of a single edge in a uniformly random linear arrangement of the vertices. It is known that \cite{Ferrer2013b}
\begin{eqnarray}
\E_{rla}[d^2] = \frac{1}{6}n(n+1)
\label{2nd_moment_of_sum_of_edge_lengths_one_edge_equation}\\
\V_{rla}[d] = \frac{1}{18}(n+1)(n-2)
\label{variance_of_sum_of_edge_lengths_one_edge_equation}
\end{eqnarray}
for $n \geq 2$.
Expressions equivalent to Eq. \ref{1st_moment_of_sum_of_edge_lengths_equation}, \ref{2nd_moment_of_sum_of_edge_lengths_one_edge_equation} and 
\ref{variance_of_sum_of_edge_lengths_one_edge_equation} were obtained in the pioneering work of Z\"ornig for a variable $d' = d - 1$ when investigating the distance between like elements in random permutations of sequences \cite{Zornig1984a}. 
Here we aim to obtain simple formulae of $\E_{rla}[D]$, $\E_{rla}[D^2]$ and $\V_{rla}[D]$ for an arbitrary network. 

The remainder of the article is organized as follows. Section \ref{preliminaries_section} {presents compact formulae for $\E_{rla}[D]$, examines the distribution of $D$ in complete graphs in detail and introduces elementary definitions and concepts. Section \ref{second_moment_section} presents compact formulae for $\E_{rla}[D^2]$ and $\V_{rla}[D]$ showing that they depend only on $n$, $m$ and $\left<k^2\right>$, the 2nd moment about zero of degree of the network under consideration, defined as
\begin{equation}
\left< k^2 \right> = \frac{1}{n} \sum_{i=1}^n k_i^2,
\label{2nd_moment_of_degree_equation}
\end{equation}
where $k_i$ is the degree of the $i$-th vertex.
% Notice that $\left<k^2\right>$ is actually an expectation, namely the expected degree of a vertex choosen uniformly at random. 
Section \ref{Erdos_Renyi_graph_section} analyses the evolution of $\V_{rla}[D]$ as $m$ increases in Erd\H{o}s-R\'enyi graphs with a constant number $n$ of vertices. A bell-shape peaking when the density of links is about 1/2 is found. 
Section  \ref{uniformly_random_labelled_trees_section} focuses on trees because of their interest for research on edge lengths \cite{Liu2017a}. That section explores the range of variation of $\V_{rla}[D]$ (delimited below by a linear tree and above by a star tree) and its linear dependence on $\left< k^2 \right>$ when $n$ is constant. It also investigates the scaling of the expected $\V_{rla}[D]$ as a function of $n$, that is asymptotically a power-law of $n$. 
Section \ref{upper_bounds_section} presents some elementary upper bounds of $D$. Section \ref{applications_section} outlines various empirical and theoretical applications of the theoretical results on the distribution of $D$ that have been obtained in the preceding sections. Finally, Section \ref{discussion_section}, reviews and discusses the findings.

\section{Preliminaries}

\label{preliminaries_section}

Here we present a simple derivation of $\E_{rla}[D]$, study the distribution of $D$ in complete graphs with minimal mathematical tools and introduce some notation and the concept of the number of independent edges.
  
\subsection{The first moment about zero}

A simple formula for $\E_{rla}[D]$, the 1st moment about zero of $D$  in a uniformly random linear arrangement, for an arbitrary network is not forthcoming to our knowledge. However, it is easy to derive. 
Thanks to Eq. \ref{sum_of_edge_lengths_equation}, one has
\begin{equation*}
\E_{rla}[D] = \sum_{i=1}^m \E_{rla}[d_i]. 
\end{equation*}
It is well known that the expected length of an arbitrary edge in a uniformly random linear arrangement is \cite{Ferrer2004b}
\begin{equation*}
\E_{rla}[d] =\frac{n+1}{3}
\end{equation*}
for $n\geq 2$.
Therefore,
\begin{equation}
\E_{rla}[D] = \frac{n+1}{3} m
\label{1st_moment_of_sum_of_edge_lengths_equation}
\end{equation}
for $n\geq 0$ ($n <2$ implies  $m = 0$ which in turn produces $\E_{rla}[D]=0$ as expected). In a tree with $n\geq 1$, one has $m = n - 1$ and then
\begin{equation}
\E_{rla}[D] = \frac{n^2 - 1}{3}
\label{1st_moment_of_sum_of_edge_lengths_tree_equation}
\end{equation}
for $n \geq 1$, a result already obtained in previous work \cite{Ferrer2016d,Ferrer2013b}.

\subsection{The distribution of $D$ in a complete graph.}

\label{complete_graph_subsection}

In general, we will use $X(G)$ to indicate a property $X$ over an arbitrary graph $G$.
Then $m(G)$ is the number of edges of graph $G$. Let $D(G)$ be the sum of dependency lengths of a linear arrangement of an arbitrary graph. Let $\completegraph$ be a complete graph. It is easy to see that 
\begin{equation}
m(\completegraph) = {n \choose 2}
\label{edges_complete_graph_equation}
\end{equation}
and that $D(\completegraph)$ is constant (it does not depend on the linear arrangement). The latter follows from the fact that there are at most $n - d$ edges of length $d$ in a linear arrangement of vertices \cite{Ferrer2013b} and a complete graph takes exactly $n - d$ edges of length $d$ for $1 \leq d \leq n - 1$, regardless of the ordering of the vertices.
Therefore, Eqs. \ref{1st_moment_of_sum_of_edge_lengths_equation} and \ref{edges_complete_graph_equation} produce
\begin{eqnarray}
D(\completegraph) & = & \E_{rla}[D(\completegraph)] \nonumber \\ 
              & = & \frac{n+1}{3} m(\completegraph) \nonumber \\
              & = & \frac{n+1}{3} {n \choose 2} \nonumber \\
              & = & \frac{1}{6}(n+1)n(n-1). \label{sum_of_edge_lengths_complete_graph_equation}
\end{eqnarray}  
The fact that $D(\completegraph)$ is constant implies that $\V_{rla}[D(\completegraph)] = 0$. As 
\begin{equation*}
\V_{rla}[D] = \E_{rla}[D^2] - \E_{rla}[D]^2,  
\end{equation*}
we have that
\begin{eqnarray}
\E_{rla}[D(\completegraph)^2] & = & \E_{rla}[D(\completegraph)]^2 \nonumber \\
                   & = & \left[ \frac{n+1}{3}{n \choose 2} \right]^2 \\
                   & = & \left[ \frac{1}{6}(n+1)n(n-1) \right]^2 \label{2nd_moment_of_sum_of_edge_lengths_complete_graph_equation}
\end{eqnarray}
for $n \geq 0$.

\subsection{The number of independent edges}

An important concept for the derivations of next sections is $q$, the number of independent pairs of edges \cite{Piazza1991a}. Two edges are said to be independent if they are not adjacent, namely, they do not share any vertex \cite[p. 4]{Bollobas1998a}. The number of pairs of different edges that can be made is \cite{Piazza1991a}
\begin{eqnarray}
Q_1 & = & {m \choose 2} \nonumber \\
        & = & \frac{1}{2}m(m-1). \label{subnumber1_of_independent_edge_pairs_equation}
\end{eqnarray}
Then $q \leq Q_1$. Besides, there are 
\begin{equation*}
{k_i \choose 2}
\end{equation*}
pairs of different edges that share vertex $i$.  
The total number of pairs of different edges that share one vertex is 
\begin{eqnarray}
Q_2 & = & \sum_{i=1}^n {k_i \choose 2} \nonumber \\
        & = &  \frac{1}{2}\left[ \sum_{i=1}^n k_i^2 - \sum_{i=1}^n k_i \right]. \label{subnumber2_of_independent_edge_pairs_equation}
\end{eqnarray}
Applying the definition of the 2nd moment about zero of degree (Eq. \ref{2nd_moment_of_degree_equation}) and the handshaking lemma \cite[p. 4]{Bollobas1998a}, namely
\begin{equation*}
2m = \sum_{i=1}^n k_i,
\end{equation*} 
one finally obtains 
\begin{equation*}
Q_2 = \frac{1}{2}n\left<k^2\right> - m.
\end{equation*}

Combining the results above (Eqs. \ref{subnumber1_of_independent_edge_pairs_equation} and \ref{subnumber2_of_independent_edge_pairs_equation}), $q$ can be defined as \cite{Piazza1991a},   
\begin{eqnarray}
q & = & Q_1 - Q_2 \nonumber \\
  & = & \frac{1}{2}\left[ m(m+1) - n \left<k^2\right> \right].
\label{number_of_independent_edge_pairs_equation}
\end{eqnarray}
%  \textcolor{red}{for $n,m \geq 0$.}
$q$ is also known as the size of the set of pairs of edges that may cross in a linear arrangement \cite{Gomez2016a}.
Hereafter we will interpret $n \left<k^2\right>$ as equivalent to 
\begin{equation}
\sum_{i=1}^n k_i^2
\end{equation}
so that the product is properly defined even when $n=0$ (calculating $n \left<k^2\right>$ from $n$ and $\left<k^2\right>$ separately is problematic because $\left<k^2\right>$ is a mean that is not defined when $n=0$; recall Eq. \ref{2nd_moment_of_degree_equation}). As a result, Eq. \ref{number_of_independent_edge_pairs_equation} is valid for $n\geq 0$.

The definition of $m(\completegraph)$ in Eq. \ref{edges_complete_graph_equation} and 
\begin{equation*}
\left<k^2\right>(\completegraph) = (n - 1)^2
\end{equation*}
transform Eq. \ref{number_of_independent_edge_pairs_equation} into
\begin{equation}
q(\completegraph) = \frac{1}{8}n(n-1)(n-2)(n-3)
\label{pairs_of_independent_edges_equation}
\end{equation}
after some algebra. 
  
\section{The second moment about zero and the variance}

\label{second_moment_section}

By definition, we have (recall Eq. \ref{sum_of_edge_lengths_equation})
\begin{eqnarray*}
\fl \E_{rla}[D^2] & = & \E\left[\left(\sum_{i=1}^m d_i \right)^2 \right] \\
\fl              & = & \E_{rla}[d_1 d_1] + \E_{rla}[d_1 d_2] +...+ \E_{rla}[d_i d_j] + ... + \E_{rla}[d_m d_{m-1}] + \E_{rla}[d_m d_m]. 
\end{eqnarray*} 
The terms $\E_{rla}[d_i d_j]$ can be classified according to $\phi$, the number of vertices shared between the $i$-th and the $j$-th vertex: $\phi = 0$ if the edges do not share any vertex, $\phi = 1$ if the edges share just one vertex and $\phi = 2$ if the edges are identical. 
This allows one to define $\E_\phi$ as the expectation of $d_i d_j$ when the $i$-th edge and the $j$-th edge share $\phi$ vertices and express the second moment about zero as 
\begin{equation}
\E_{rla}[D^2] = \sum_{\phi = 0}^2 f_\phi \E_\phi,
\label{raw_2nd_moment_of_sum_of_edge_lengths_equation}
\end{equation} 
where $f_\phi$ is the number of terms of type $\phi$. 
Obviously, 
\begin{equation}
f_2 = m.
\label{raw_frequency_type_2_equation}
\end{equation}
It is easy to see that 
\begin{equation}
f_0 = 2q.
\label{raw_frequency_type_0_equation}
\end{equation}
The fact that 
\begin{equation}
m^2 = \sum_{\phi = 0}^2 f_\phi
\label{total_frequency_equation}
\end{equation}
combined with $f_2 = m$ and $f_0 = 2q$ gives
\begin{equation}
f_1 = m(m-1) - 2q.
\label{raw_frequency_type_1_equation}
\end{equation}

Recalling the formula for $m(\completegraph)$ (Eq. \ref{edges_complete_graph_equation}) and that of $q(\completegraph)$ (Eq. \ref{pairs_of_independent_edges_equation}), it is easy to see that 
\begin{eqnarray}
f_0(\completegraph) & = & \frac{1}{4}n(n-1)(n-2)(n-3) \label{frequency_type_0_complete_graph_equation} \\
f_1(\completegraph) & = & n(n-1)(n-2) \label{frequency_type_1_complete_graph_equation} \\ 
f_2(\completegraph) & = & \frac{1}{2}n(n-1) \label{frequency_type_2_complete_graph_equation}
\end{eqnarray}
after some algebra.

Now we turn our attention to the calculation of $\E_0$, $\E_1$ and $\E_2$ for an arbitrary graph. The calculation of $\E_2$ is straightforward. Note that $\E_2 = \E_{rla}[d^2]$ where $d$ is the length of an arbitrary edge. Eq. \ref{2nd_moment_of_sum_of_edge_lengths_one_edge_equation} gives 
\begin{equation}
\E_2 = \frac{1}{6}n(n+1)
\label{expectation_type_2_equation}
\end{equation}
for $n \geq 2$.

The calculation of $\E_1$ is more elaborate and requires enumerating all the possible linear arrangements of the vertices of the $i$-th and the $j$-th edge when they share one vertex. 
The number of linear arrangements where 
\begin{enumerate}
\item
The shared vertex is located in between the two vertices
\item
The $i$-th edge appears first
\end{enumerate} 
is
\begin{eqnarray*}
A & = & \sum_{d_i=1}^{n-2} \sum_{d_j=1}^{n-1-d_i} (n - d_i - d_j) \\
  & = & \frac{1}{6} n(n-1)(n-2).
\end{eqnarray*} 
The total number of linear arrangements where the shared vertex is located in between the two vertices is thus $2A$. 

The number of linear arrangements where the shared vertex is located after the other two vertices and the other vertex of the $i$-th edge appears first is
\begin{eqnarray*}
B & = & \sum_{d_i=2}^{n-1} \sum_{d_j=1}^{d_i - 1} (n - d_j) \\
  & = & \frac{1}{6} n(n-1)(n-2).
\end{eqnarray*} 
The number of linear arrangements where the shared vertex is located either after or before the other two vertices is thus $4B$.
Therefore we conclude that the total number of linear arrangements of the vertices of two edges that share one vertex is 
\begin{equation*}
T = 2A+4B.
\end{equation*}
As $\E_1$ is the average value of $d_i d_j$ over all linear arrangements of the three vertices, we have that 
\begin{eqnarray}
\E_1 & = & \frac{2}{T}(A'+2B') \nonumber \\
             & = &  \frac{A'+2B'}{A+2B}, \label{raw_expectation1_equation}
\end{eqnarray}
where
\begin{eqnarray*}
A' & = & \sum_{d_i=1}^{n-2} \sum_{d_j=1}^{n-1-d_i} (n - d_i - d_j) d_i d_j\\
   & = & \frac{1}{120} (n+2)(n+1)n(n-1)(n-2)
\end{eqnarray*}
and 
\begin{eqnarray*}
B' & = & \sum_{d_i=2}^{n-1} \sum_{d_j=1}^{d_i - 1} (n - d_i) d_i d_j \\
   & = & \frac{1}{120} (3n+1)(n+1)n(n-1)(n-2).
\end{eqnarray*} 
Applying the expressions for $A$, $B$, $A'$ and $B'$ that have been obtained above to Eq. \ref{raw_expectation1_equation}, one obtains
\begin{equation}
\E_1 = \frac{1}{60}(n+1)(7n+4)
\label{expectation_type_1_equation}
\end{equation} 
for $n\geq$ 3.

To calculate $\E_0$, we take the definition of $\E_{rla}[D^2]$ in Eq. \ref{raw_2nd_moment_of_sum_of_edge_lengths_equation} and obtain 
\begin{equation}
\E_0 = \frac{\E_{rla}[D^2] - f_1 \E_1 - f_2 \E_2}{f_0}.
\label{precursor_equation}
\end{equation}
Recall that $\E_\phi$ is simply the expected value of the product of two lengths from a couple of edges that share $\phi$ vertices. The only constraint on $\E_\phi$ is that $n \geq 4 - \phi$.
Notice that, given an $n$ that satisfies such a constraint, $\E_\phi$ is defined independently from the kind of graph under consideration. % This is the rationale of the computional validation of the formulae for the $\E_\phi$'s (\ref{validation_section}) and their alternative derivation  procedure showing that their values depend only on $n$ (\ref{alternative_derivation_section}). 
  % Thus, the $\E_\phi$'s have a constant value in every graph where the corresponding type exists, namely $f_\phi >0$ (for instance, type 0 is missing in a star tree because $f_0  = q = 0$ \cite{Gomez2016a}).
Then, we will derive $\E_0$ borrowing $\E_{rla}[D^2]$, $f_1$ and $f_2$ from a complete graph obtaining a value of $\E_0$ that is valid for arbitrary graphs.
Applying the values of $\E_{rla}[D(\completegraph)^2]$ (Eq. \ref{2nd_moment_of_sum_of_edge_lengths_complete_graph_equation}), $f_0(\completegraph)$ (Eq. \ref{frequency_type_0_complete_graph_equation}), $f_1(\completegraph)$ (Eq. \ref{frequency_type_1_complete_graph_equation}) and $f_2(\completegraph)$ (Eq. \ref{frequency_type_2_complete_graph_equation}) as well as the values of $\E_1$ (Eq. \ref{expectation_type_1_equation}) and $\E_2$ (Eq. \ref{expectation_type_2_equation}), one gets
\begin{equation}
\E_0 = \frac{1}{45}(n+1)(5n+4)
\label{expectation_type_0_equation}
\end{equation}
for $n\geq 4$ after some algebra.

Now we aim to find a compact formula for $\E_{rla}[D^2]$. The definitions of $\E_\phi$ (Eqs. \ref{expectation_type_0_equation}, \ref{expectation_type_1_equation} and \ref{expectation_type_2_equation}) transform Eq. \ref{raw_2nd_moment_of_sum_of_edge_lengths_equation} into
\begin{equation*}
\E_{rla}[D^2] = \frac{n+1}{3} \left[ \frac{1}{15}(5n+4)f_0 + \frac{1}{20}(7n+4)f_1 + \frac{1}{2}nf_2 \right].
\end{equation*}
The fact that $f_1 = m^2 - f_0 - f_2$ (recall Eq. \ref{total_frequency_equation}), gives  
\begin{equation*}
\E_{rla}[D^2] = \frac{n+1}{3} \left[ \frac{1}{60}(4-n)f_0 + \frac{1}{20}(3n- 4)f_2 + \frac{1}{20}(7n +4) \right].
\end{equation*}
Applying $f_2 = m$ and $f_0 = 2q = m(m+1) - n\left<k^2 \right>$ (recall Eq  \ref{number_of_independent_edge_pairs_equation}), one obtains 
\begin{equation}
\E_{rla}[D^2] = \frac{n+1}{45} \left[ m(m(5n + 4) + 2(n-1)) + \left(\frac{n}{4} - 1 \right)n \left<k^2 \right> \right].
\label{2nd_moment_of_sum_of_edge_lengths_equation}
\end{equation}  
% \textcolor{red}{for $n \geq 0$.}
after some work. % Eq. \ref{2nd_moment_of_sum_of_edge_lengths_equation} shows that $\E_{rla}[D^2]$ depends only on $n$, $m$ and $\left<k^2\right>$.
The variance of $D$ is
\begin{eqnarray}
\V_{rla}[D] & = & \E_{rla}[D^2] - \E_{rla}[D]^2 \nonumber \\
     & = & \E_{rla}[D^2] - \frac{1}{9}m^2(n+1)^2 \nonumber \\
     & = & \frac{n+1}{45} \left[ m(2(n-1) - m ) + \left(\frac{n}{4} - 1 \right)n \left<k^2 \right> \right]. \label{variance_of_sum_of_edge_lengths_equation}
\end{eqnarray}
% \textcolor{red}{for $n \geq 0$.}
Eqs. \ref{2nd_moment_of_sum_of_edge_lengths_equation} and \ref{variance_of_sum_of_edge_lengths_equation} show that $\E_{rla}[D^2]$ and $\V_{rla}[D]$ depend only on $n$, $m$ and $\left<k^2 \right>$.
It is easy to see that Eq. \ref{2nd_moment_of_sum_of_edge_lengths_equation} and \ref{variance_of_sum_of_edge_lengths_equation} are valid for $n,m \geq 0$ because 
\begin{itemize}
\item
The $f_\phi$'s are valid for $n,m \geq 0$.
\item
$E_\phi$ is valid only for $n \geq 4 - \phi$ but  $f_\phi = 0$ for $n < 4 - \phi$.
\item
In Eq \ref{raw_2nd_moment_of_sum_of_edge_lengths_equation}, the product by $f_\phi$ warrants that an invalid value of $\E_\phi$ will have zero contribution. 
\end{itemize}
In a tree with $n \geq 1$,  $m = n - 1$ and then
\begin{eqnarray}
\E_{rla}[D^2] & = & \frac{n+1}{45} \left[ (n-1)^2 (5n+6) + \left(\frac{n}{4} - 1 \right)n \left<k^2 \right> \right] \label{2nd_moment_of_sum_of_edge_lengths_tree_equation}\\
\V_{rla}[D] & = & \frac{n+1}{45} \left[ (n-1)^2 + \left(\frac{n}{4} - 1 \right)n \left<k^2 \right> \right].
\label{variance_of_sum_of_edge_lengths_tree_equation}
\end{eqnarray}
for $n \geq 1$.
Thus, $\E_{rla}[D^2]$ and $\V_{rla}[D]$ are completely determined by $n$ and $\left<k^2 \right>$ in trees.
See \ref{validation_section} for the procedure that we have used to check the theoretical results obtained so far and \ref{alternative_derivation_section} for an alternative derivation of the variance of $D$.

Let us calculate $\V_{rla}[D]$ for the network in Fig. \ref{syntactic_dependency_tree_figure}. Eq. \ref{2nd_moment_of_degree_equation} can be expressed equivalently as 
\begin{equation}
\left< k^2 \right> = \frac{1}{n} \sum_{k=1}^{n-1} n(k) k^2,
\end{equation}
where $n(k)$ is the number of vertices of degree $k$. The summary of vertex degrees in Table \ref{degrees_table} yields $\left< k^2 \right> = 88/17$. Applying this result and $n = 17$ to Eq. \ref{variance_of_sum_of_edge_lengths_tree_equation} one obtains $\V_{rla}[D] = 1084/5$. 
Table \ref{syntactic_dependency_tree_table} summarizes the statistical properties of the network.

\begin{table}
\caption{\label{degrees_table} $n(k)$, the number of vertices of degree $k$, for the network in Fig. \ref{syntactic_dependency_tree_figure}. $n(k)=0$ for $k > 5$. }
\begin{indented}
\item[]
\begin{tabular}{@{}rr}
\br
$k$ & $n(k)$ \\ % \hline
\mr
1 & 9 \\
2 & 6 \\
3 & 0 \\
4 & 0 \\
5 & 2 \\
\br
\end{tabular}
\end{indented}
\end{table}
 
\section{Erd\H{o}s-R\'enyi graphs}

\label{Erdos_Renyi_graph_section}

Let us consider ${\cal G}_{n,m}$, the ensemble of graphs of $n$ vertices and $m$ edges where all the 
\begin{equation*}
{{n \choose 2} \choose m}
\end{equation*} 
distinct graphs of $m$ edges are equally likely \cite{Erdos1959a}. A sibling ensemble is ${\cal G}_{n,\pi}$, that was introduced by Gilbert and that consists of graphs where a pair of different vertices are linked with probability $\pi$ independently from other vertex pairs \cite{Bollobas2002a}. We call these two ensembles siblings because they behave similarly and are almost interchangeable when $m \approx \pi n$ \cite{Bollobas2002a}. 

We aim to predict the expected value of $\E_{rla}[D]$, $\E_{rla}[D^2]$, $\V_{rla}[D]$ in ${\cal G}_{n,m}$ theoretically. Suppose that $\E_{n,m}$ is the expectation operator over the ensemble ${\cal G}_{n,m}$. 
Then Eq. \ref{1st_moment_of_sum_of_edge_lengths_equation} gives 
\begin{eqnarray*}
\E_{n,m}\left[\E_{rla}[D]\right] & = & \E_{n,m}\left[\frac{n+1}{3}m \right] \\
                                       & = & \frac{n+1}{3}m
\end{eqnarray*}
trivially since both $n$ and $m$ are constant.
 
Let us consider the case of $\E_{n,m}\left[\V_{rla}[D]\right]$. Then Eq. \ref{variance_of_sum_of_edge_lengths_equation} gives
\begin{equation}
\E_{n,m}\left[\V_{rla}[D]\right] = \frac{n+1}{45} \left[ m(2(n-1) - m) + \left(\frac{n}{4} - 1 \right)n \E_{n,m}\left[\left<k^2 \right>\right] \right].
\label{raw_expected_variance_Erdos_Renyi_equation}
\end{equation}
One may calculate $\E_{n,m}[\left<k^2 \right>]$ knowing that the probability that a vertex has degree $k$ in a graph from ${\cal G}_{n,m}$ is \cite[bottom of p. 58]{Erdos1960a}
\begin{equation*}
p(k) = \frac{{n-1 \choose k}{{n-1 \choose 2}\choose m - k}}{{{n \choose 2}\choose m}}.
\end{equation*}
Notice that $p(k) = 0$ if $k > m$ as expected.
Then $\E_{n,m}[\left<k^2 \right>]$ can be replaced by    
\begin{equation*}
\left<k^2\right> = \sum_{k=1}^{n-1}p(k) k^2
\end{equation*}
as all graphs of the ensemble are equally likely. 
Unfortunately, a closed form formula is not available to our knowledge.

Figure \ref{Erdos_Renyi_graph_figure} compares the predicted $\V_{rla}[D]$ against estimates via a Monte Carlo procedure for various values of $n$ and all possible values of $m$. To ease visualization, $D$ and $m$ are normalized dividing them by their respective maximum value, that are achieved by a complete graph.  
Dividing $m$ by its maximum value ($m(\completegraph)$) one obtains $\delta$, the density of links (Eq. \ref{density_of_links_equation}).  
$\V_{rla}[D]$ is normalized dividing it by $D(\completegraph)$. It can be seen that the theoretical curve matches the simulations accurately. 
Figure \ref{Erdos_Renyi_graph_figure} shows that $\V_{rla}[D]$ is minimized by extreme values of $m$ as expected (in an empty graph or a complete graph, $\V_{rla}[D] = 0$) and that $\V_{rla}[D]$ is maximized by intermediate values of $m$. 

For each value of $n$ and $m$, the Monte Carlo procedure estimates $\V_{rla}[D]$ using  an unbiased estimator over $T=10^4$ replicas. A naive procedure may estimate $\V_{rla}[D]$ generating $T$ random graphs for every possible value of $m$. Instead, we use a more efficient procedure:  
\begin{enumerate}
\item
Generate a vector containing the ${n \choose 2}$ edges of the complete graph.
\item
Repeat $T$ times the following subprocedure
   \begin{enumerate}
   \item
   Shuffle the content of the vector producing a uniformly random permutation of its content. 
   \item
   Generate a random graph of $m$ vertices by choosing one edge at a time. The $m$ first edges of the vector define a random graph of $m$ edges. 
   \end{enumerate}
\item
Calculate the expected $\V_{rla}[D]$ as a function of $m$ as the mean value of $\V_{rla}[D]$ over the $T$ random graphs that have been produced for each value of $m$ in the previous step. 
\end{enumerate}

\begin{figure}
\begin{center}
\includegraphics[scale = 0.6]{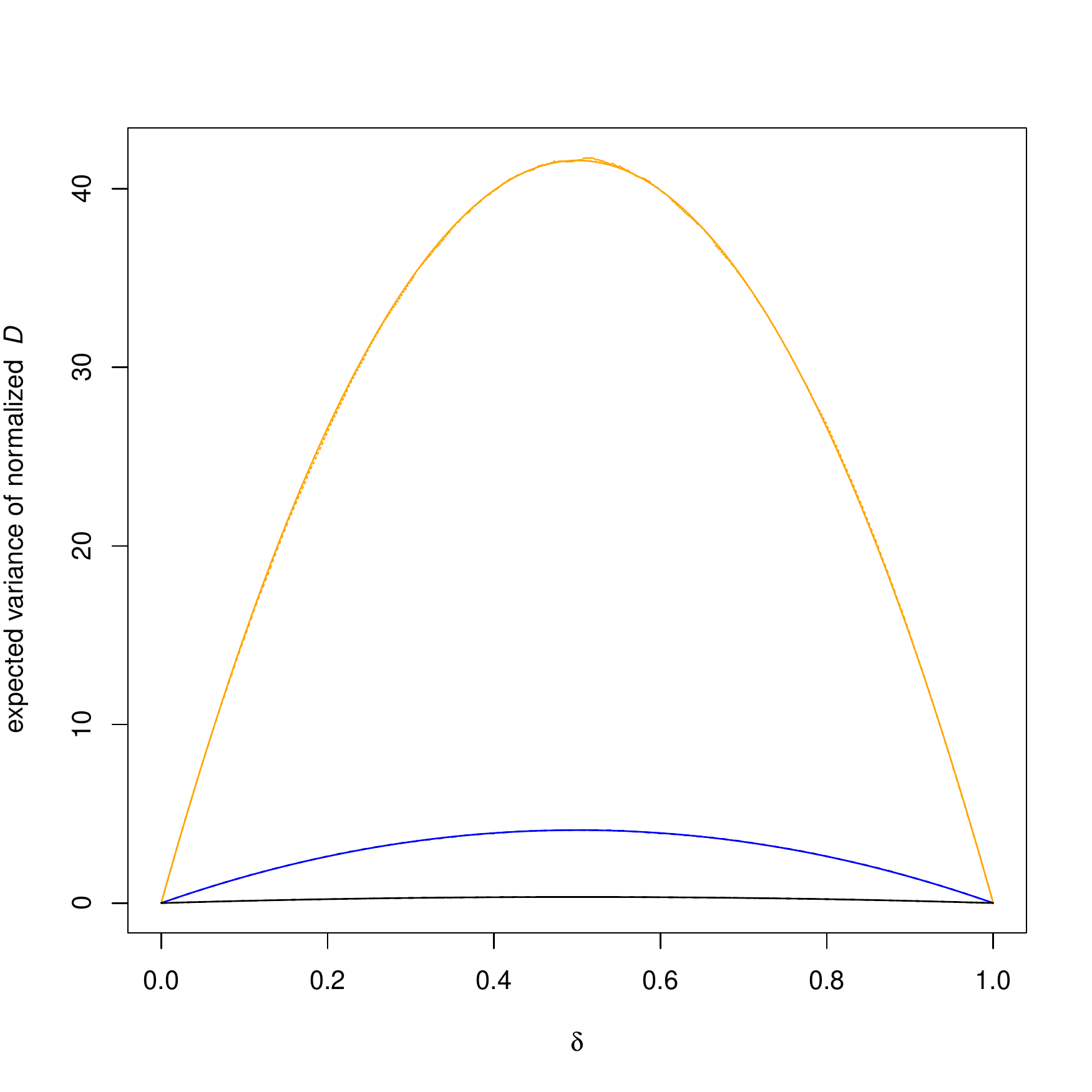}
\end{center}
\caption{\label{Erdos_Renyi_graph_figure} The expected variance of normalized $D$ in ${\cal G}_{n,m}$ ($\E_{n,m}[\V_{rla}[D]/D(\completegraph)]$) as a function of the density of links ($\delta$) for different values of $n$ ($n=10$ in black, $n=100$ in blue and $n=1000$ in orange). Solid lines are used for theoretical predictions, dashed lines are used for estimates via simulation. Dashed lines can hardly be seen due to the accuracy of the theoretical predictions and the high number of replicas used. The values of $\delta$ are obtained applying $\delta = m/m(\completegraph)$ to all possible values of $m \in [1, m(\completegraph)]$. }
\end{figure}

We consider alternative ways to approximate $\E_{n,m}[\left<k^2 \right>]$ to have a compact (albeit approximated) formula for $\left<k^2\right>$ and to be able to locate the maxima of $\V_{rla}[D]$ for intermediate values of $m$ precisely (Fig. \ref{Erdos_Renyi_graph_figure}).
The close relationship between ${\cal G}_{n,m}$ and ${\cal G}_{n,\pi}$ \cite{Bollobas2002a} suggests that 
\begin{equation*}
\E_{n,m}\left[\left<k^2 \right>\right] \approx \E_{n,\pi}\left[\left<k^2 \right>\right]
\end{equation*}
with $\pi = \delta$, 
where 
\begin{equation}
\delta = \frac{m}{{n \choose 2}}
\label{density_of_links_equation}
\end{equation}
is the density of links. This approximation is convenient because $\left<k^2\right>$ has a simple closed form formula in ${\cal G}_{n,\pi}$. In particular, $k$ follows a binomial distribution with parameters $n - 1$ and $\pi$ in a graph from ${\cal G}_{n,\pi}$, namely 
\begin{equation*}
p(k) = {n - 1 \choose k} \pi^k (1-\pi)^{n - 1 - k}
\end{equation*}
and then
\begin{eqnarray*}
\left<k\right> & = & (n - 1)\pi \\
\V[k]                & = & (n - 1)\pi(1 - \pi) \\
\left<k^2\right> & = & \V[k] + \left<k\right>^2 \\
                 & = & (n-1)\pi[(n-2)\pi + 1].
\end{eqnarray*}
This is why we call it a binomial degree approximation. 

Choosing the most likely value of $\pi$ for a graph in ${\cal G}_{n,m}$, namely, $\pi = \delta$ one obtains
\begin{eqnarray}
\E_{n,m}\left[\left<k^2 \right>\right] & \approx & (n-1)\delta((n-2)\delta + 1) \nonumber \\
                           & = & (n-1)\frac{m}{{n \choose 2}}\left((n-2)\frac{m}{{n \choose 2}} + 1 \right) \label{second_moment_of_degree_approximation_equation}\\
                           & = & \frac{4(n-2)}{(n - 1)n^2}m^2 - \frac{2m}{n}. \nonumber 
\end{eqnarray}
Applying the last result to Eq. \ref{raw_expected_variance_Erdos_Renyi_equation}, one gets 
\begin{equation}
\E_{n,m}\left[\V_{rla}[D]\right] \approx \frac{(n+1)m}{45} \left[\frac{8 - 5n}{n(n-1)}m + 2\left(\frac{5n}{4} -2 \right) \right].
\label{expected_variance_Erdos_Renyi_equation}
\end{equation}
Equating 
\begin{equation*}
\frac{d \E_{n,m}\left[\V_{rla}[D]\right]}{dm} \approx 2\frac{n+1}{45} \left[\frac{8 - 5n}{n(n-1)}m + \frac{5n}{4} -2 \right]
% \label{1st_derivative_of_expected_variance_Erdos_Renyi_equation}
\end{equation*}
to zero, one finds that $\E_{n,m}[\V_{rla}[D]]$ has a critical point approximately at  
\begin{equation*}
m^* =  \frac{n(n-1)}{4}.
\end{equation*}
The second derivative is negative, indicating that the critical point is a maximum. 
Although $m^*$ has been obtained via a binomial degree approximation for an Erd\H{o}s-R\'enyi, Fig. \ref{Erdos_Renyi_graph_approximation_figure} shows that the approximation is accurate.

% As $m$ is a natural number by definition, one has that the maximum is unique when $m^*$ is integer. When $m^*$ is not integer, there are two maxima, one at $\left\lfloor m^* \right\rfloor$ and the other at $\left\lceil m^* \right\rceil$.
% Therefore, $\V_{rla}[D]$ is maximized when $\delta \approx 1/2$. 
As $m$ is a natural number by definition, one has that the maximum is unique when $m^*$ is integer. When $m^*$ is not integer, the maxima could be located at $\left\lfloor m^* \right\rfloor$, $\left\lceil m^* \right\rceil$ or both. As $m^*$ has been obtained via an approximation, we are conservative and conclude that $\V_{rla}[D]$ is maximized when $\delta \approx 1/2$.

We already know that $\V_{rla}[D]$ reaches two global minima when $m$ is minimum, i.e. $\delta = 0$, and also when $m$ is maximum, i.e. $\delta = 1$.

In order to derive $\E_{n,m}\left[\E_{rla}[D^2]\right]$, we could apply the method above to Eq. \ref{2nd_moment_of_sum_of_edge_lengths_equation}.
% wolfram: ((n+1)/45) [ m(m(5n + 4) + 2(n-1)) + (n/4 - 1 )n [{4(n-2)}/{(n - 1)n^2} m^2 - 2m/n] ] 
Instead, we choose a faster track to obtain a compact formula, based on the fact that
\begin{equation*}
\E_{n,m}\left[\E_{rla}[D^2]\right]  =  \E_{n,m}\left[\V_{rla}[D]\right] - \E_{n,m}\left[\E_{rla}[D^2]\right].
\end{equation*}
Applying Eqs. \ref{1st_moment_of_sum_of_edge_lengths_equation} and \ref{expected_variance_Erdos_Renyi_equation}, one obtains
\begin{equation}
% ((n+1)m/45)[(8-5n)m/(n(n-1))+2(5n/4 - 2)] + [(n+1)m/3]^2
\E_{n,m}\left[\E_{rla}[D^2]\right] \approx \frac{m(n+1)}{90}\left[\frac{2m (5n(n^2-2)+8)}{n(n-1)} + 5n - 8 \right]
\end{equation}
after some routine calculations. It is easy to see that $\E_{n,m}\left[\E_{rla}[D^2]\right]$ is a monotonically increasing function of $m$ when $n$ is kept constant, contrary to the bell-shape behavior of $\E_{n,m}\left[\V_{rla}[D]\right]$.

\begin{figure}
\begin{center}
\includegraphics[scale = 0.6]{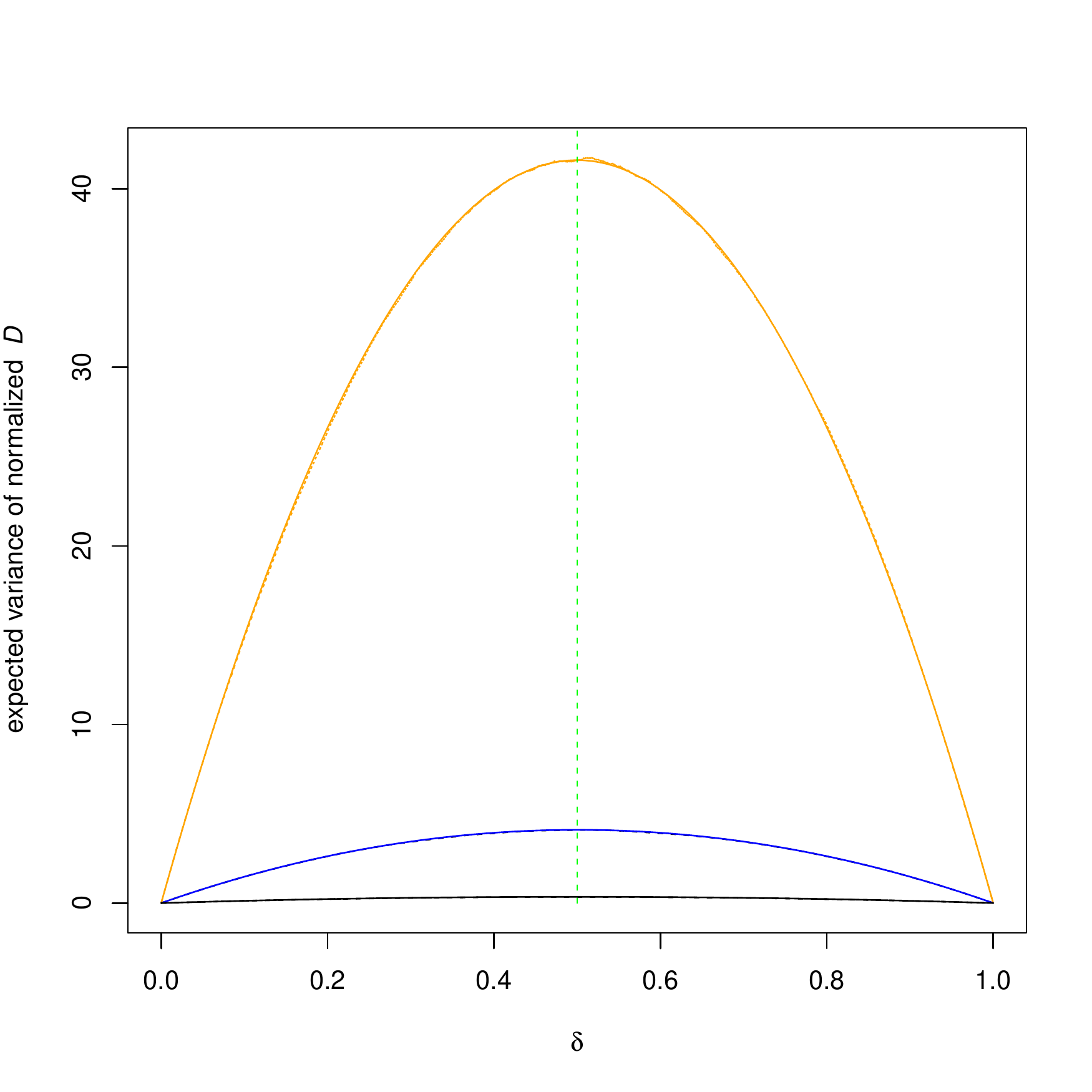}
\end{center}
\caption{\label{Erdos_Renyi_graph_approximation_figure} The same as Fig. \ref{Erdos_Renyi_graph_figure} but theoretical predictions (solid lines) come from assuming a binomial degree distribution. Dashed lines can hardly be seen due to the accuracy of the theoretical predictions and the high number of replicas used. An additional vertical line (green dashed) is used to indicate the location of the maximum $\V[k]$ at $\delta \approx 1/2$. }
\end{figure}

The fact that $k$ follows approximately a Poisson distribution, i.e. \cite{Erdos1960a}
\begin{equation*}
p(k) = e^{-\lambda} \frac{\lambda^k}{k!}
\end{equation*} 
with $\lambda = 2m/n$, allows to approximate the second moment about zero of degree as 
\begin{eqnarray}
\left<k^2\right> & = & \V[k] + \left<k\right>^2 \nonumber \\
                 & \approx & \lambda(1 + \lambda) \nonumber \\
                 & = & \frac{2m}{n}\left(1 + \frac{2m}{n} \right). \label{second_moment_of_degree_Poisson_approximation_equation}
\end{eqnarray}
However, this Poisson distribution approximation is poor. One reason to suspect this is true is the large difference between Eq. \ref{second_moment_of_degree_approximation_equation} and Eq. \ref{second_moment_of_degree_Poisson_approximation_equation}, that is
\begin{equation*}
\frac{4m}{n}(m - 1).
\end{equation*}
A deeper reason is the good approximation provided by the binomial degree distribution of ${{\cal G}_{n,\pi}}$ for $\left< k^2 \right>$ in ${{\cal G}_{n,m}}$ (Fig. \ref{Erdos_Renyi_graph_approximation_figure}). It is well known that the Poisson distribution gives only a good approximation to a binomial distribution when $n$ is large and 
\begin{equation*}
\pi = m/{n \choose 2}
\end{equation*}
is small. Unfortunately, we are exploring the whole range of variation of $\pi$ as a result of our exhaustive exploration of the values of $m$.  

\section{Uniformly random labelled trees}

\label{uniformly_random_labelled_trees_section}

\begin{figure}
\begin{center}
\includegraphics[scale = 0.6]{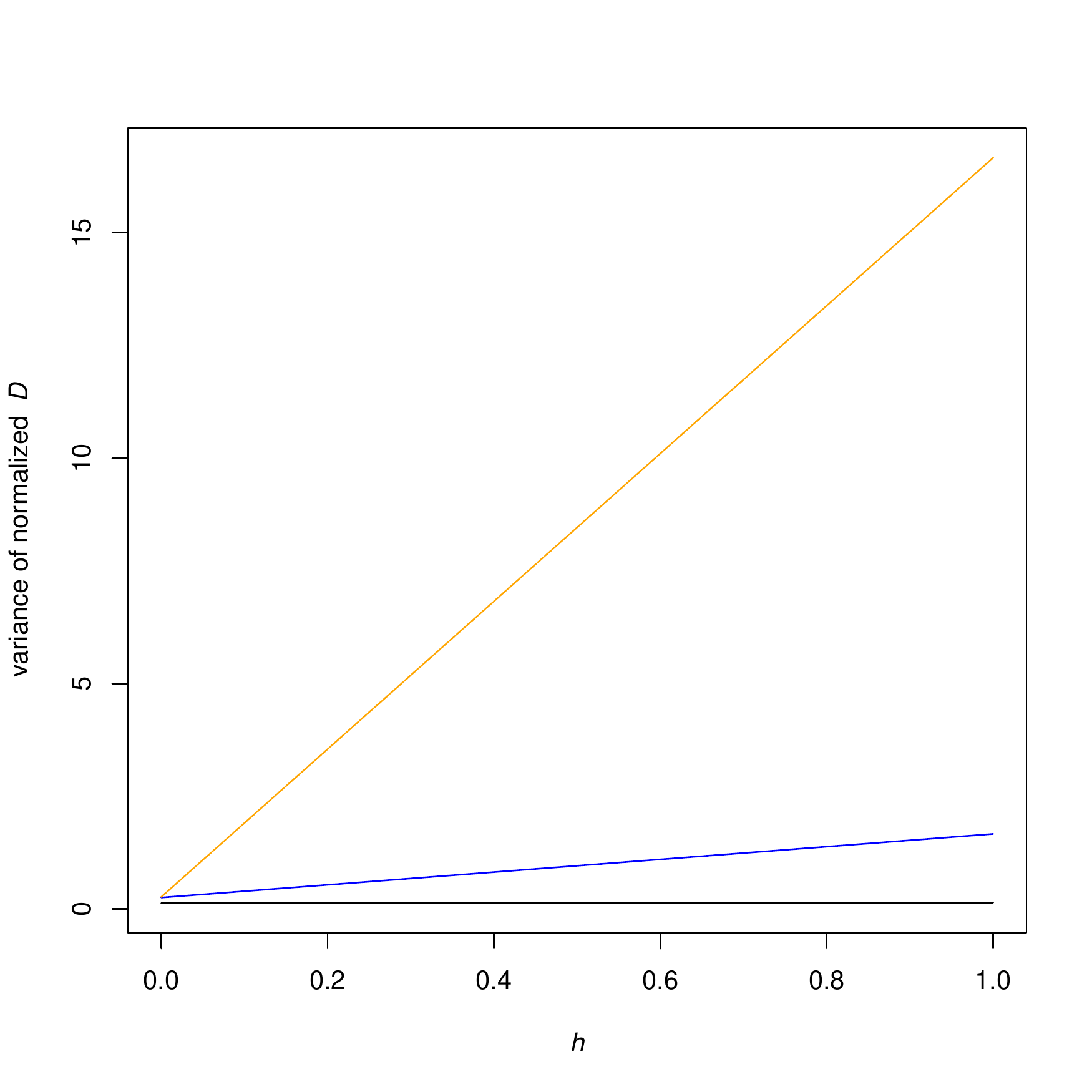}
\end{center}
\caption{\label{tree_figure} The normalized variance of $D$ ($\V_{rla}[D]/D(\completegraph)$) as a function of the hubiness coefficient of trees ($h$) for different values of $n$ ($n=5$ in black, $n=50$ in blue and $n=500$ in orange). For simplicity, the values of $h$ are obtained applying the definition of $h$ (Eq. \ref{hubiness_coefficient_equation}) to all possibles values of $\left< k^2 \right>$ under the assumption that $n\left< k^2 \right>$ can take any natural number within the interval $[n\left< k^2 \right>(\lineartree), n\left< k^2 \right>(\startree)]$. }
\end{figure}

We pay further attention to the particular case of trees given their interest for research on edge  lengths, e.g., \cite{Ferrer2004b,Liu2007a,Liu2008a,Futrell2015a,Ouyang2017a,Lei2018a}. 
If that case, $m$ becomes constant and then $\V_{rla}[D^2]$ becomes an increasing linear function of $\left<k^2\right>$ when $n$ is also constant (recall \ref{variance_of_sum_of_edge_lengths_tree_equation})). In trees where $n$ is given, it has been shown that \cite{Esteban2016a}
\begin{equation*}
\left<k^2\right>(\lineartree) \leq \left<k^2\right> \leq \left<k^2\right>(\startree),
\end{equation*}
where 
\begin{equation}
\left<k^2\right>(\lineartree) = 4 - \frac{6}{n}
\label{2nd_moment_of_degree_linear_tree}
\end{equation}
is the 2nd moment about zero of degree of $\lineartree$, a linear tree of $n$ vertices, when $n\geq 2$, and 
\begin{equation*}
\left<k^2\right>(\startree) = n - 1
\end{equation*}
is the 2nd moment about zero of degree of $\startree$, a star tree of $n$ vertices, when $n \geq 1$.
Therefore $\V_{rla}[D]$ is minimized by linear trees and maximized by star trees, i.e.
\begin{eqnarray*}
\E_{rla}[D^2](\lineartree) \leq \E_{rla}[D^2] \leq \E_{rla}[D^2](\startree) \\
\V_{rla}[D](\lineartree) \leq \V_{rla}[D] \leq \V_{rla}[D](\startree).
\end{eqnarray*}
% $\E_{rla}[D^2](\startree)$ and $\V_{rla}[D](\startree)$ can be found in Eqs. \ref{2nd_moment_of_sum_of_edge_lengths_star_tree_equation} and \ref{variance_of_sum_of_edge_lengths_star_tree_equation}. 
Applying Eq. \ref{2nd_moment_of_degree_linear_tree} to Eqs. \ref{2nd_moment_of_sum_of_edge_lengths_tree_equation} and \ref{variance_of_sum_of_edge_lengths_tree_equation} one obtains
\begin{eqnarray*}
\E_{rla}[D^2](\lineartree) = \frac{1}{90} (n + 1) (10n^3 - 6n^2 - 25n +24) \\
\V_{rla}[D](\lineartree) = \frac{1}{90} (n + 1) (n-2) (4n - 7) 
\end{eqnarray*}
for $n\geq 2$ after some algebra. Equivalent results for $\E_{rla}[D^2](\startree)$ and $\V_{rla}[D](\startree)$  are derived in \ref{validation_section} (Eqs. \ref{2nd_moment_of_sum_of_edge_lengths_star_tree_equation}
\ref{variance_of_sum_of_edge_lengths_star_tree_equation}).

We aim to explore the actual dependency between $\V_{rla}[D]$ and $\left<k^2 \right>$ in trees. To ease comparison, we will rescale these two variables. First, $D$ is normalized dividing it by $D(\completegraph)$.
Second, $\left<k^2\right>$ is normalized with the help of $h$, the hubiness coefficient, that is defined as
\begin{equation}
h = \frac{\left<k^2\right> - \left<k^2\right>(\lineartree)}{\left<k^2\right>(\lineartree) - \left<k^2\right>(\startree)}.
\label{hubiness_coefficient_equation}
\end{equation}
It is easy to see that $0 \leq h \leq 1$ ($h = 0$ in a linear tree and $h = 1$ in a star tree).
Figure \ref{tree_figure} shows the expected monotonic growth of $\V_{rla}[D]$ predicted by Eq. \ref{variance_of_sum_of_edge_lengths_equation}.

Research on the linear arrangement of trees has considered different statistical frameworks for random trees \cite{Ferrer2006d,Liu2008a,Ferrer2014c,Esteban2016a}. For simplicity, here we focus on ${\cal T}_n$, the ensemble of random trees where all the possible labelled trees of $n$ vertices are equally likely \cite{Aldous1990a,Broder1989a,Esteban2016a}. It is well known that there are $n^{n-2}$ labelled trees of $n$ vertices \cite{Cayley1889a}. Refs. \cite{Ferrer2014c,Esteban2016a} are based on this ensemble; Refs. \cite{Ferrer2006d,Liu2008a} are not. An advantage of this ensemble is the availability of results that allow one to predict the expected value of $\V_{rla}[D]$ or $\E_{rla}[D^2]$ in ${\cal T}_{n}$ theoretically.

Suppose that $\E_{rlt}$ is the expectation operator over the ensemble ${\cal T}_n$ and 
let us consider $\E_{rlt}\left[\E_{rla}[D]\right]$. 
Then Eq. \ref{1st_moment_of_sum_of_edge_lengths_tree_equation} with $m=n-1$ gives 
\begin{equation*}
\E_{rlt}\left[\E_{rla}[D]\right] = \frac{1}{3}(n+1)(n-1) 
\end{equation*}
trivially since $n$ is constant.

Now let us consider the expectation of $\V_{rla}[D]$ in that ensemble. 
Then Eq. \ref{variance_of_sum_of_edge_lengths_tree_equation} gives
\begin{equation}
\E_{rlt}\left[\V_{rla}[D]\right] = \frac{n+1}{45} \left[ (n-1)^2 + \left(\frac{n}{4} - 1 \right)n \E_{rlt}\left[\left<k^2 \right>\right] \right].
\label{raw_expected_variance_uniformly_random_tree_equation}
\end{equation}
Knowing that \cite{Moon1970a,Noy1998a,Ferrer2016d} 
\begin{equation}
\E_{rlt}\left[\left<k^2 \right>\right] = \left(1 - \frac{1}{n}\right)\left(5 - \frac{6}{n}\right),
\end{equation}
Eq. \ref{raw_expected_variance_uniformly_random_tree_equation} becomes
\begin{equation}
\E_{rlt}\left[\V_{rla}[D]\right] = \frac{(n+1)(n-1)(13n^2 - 54n + 48)}{360n}
\label{expected_variance_uniformly_random_tree_equation}
\end{equation}
after some routine calculations. Applying the same methodology to Eq. \ref{2nd_moment_of_sum_of_edge_lengths_tree_equation} one obtains
\begin{equation*}
\E_{rlt}\left[\E_{rla}[D^2]\right] = \frac{(n+2)(n+1)(n-1)(4n - 3)(5n - 4)}{180n}.
\end{equation*}
Figure \ref{scaling_of_variance_random_trees_equation} shows that the growth of $\E_{rlt}\left[\V_{rla}[D]\right]$ as a function of $n$ according to Eq. \ref{expected_variance_uniformly_random_tree_equation} matches perfectly the estimates from Monte Carlo procedure. It also shows that $\E_{rlt}\left[\V_{rla}[D]\right]$ is closer to the reference $\V_{rla}[D]$ of a linear tree than to that of a star tree.

The Monte Carlo procedure consists of estimating $\E_{rlt}\left[\V_{rla}[D]\right]$ as the sample variance of $D$ in $T = 10^6$ uniformly random trees. Each random tree is produced generating a uniformly random Pr{\"u}fer code and transforming it into a labelled tree \cite{Pruefer1918a}. Such a procedure turns out to be computationally optimal to generate a tree with given $n$ \cite[Chapter 3.3]{Alonso1995a}. A uniformly random linear arrangement of the vertices is assigned to each tree. The same method is used to estimate  $\V_{rla}[D]$ in linear trees and star trees. The only difference is that the tree is given, not generated at random.

It is easy to see from the equations above that asymptotically 
\begin{eqnarray*}
\E_{rlt}\left[\V_{rla}[D]\right], \V_{rla}[D](\lineartree) \sim n^3 \\
\E_{rla}[D^2](\startree), \V_{rla}[D](\startree), \E_{rlt}\left[\E_{rla}[D^2]\right], \E_{rla}[D^2](\lineartree) \sim n^4.
\end{eqnarray*}

\begin{figure}
\begin{center}
\includegraphics[scale = 0.6]{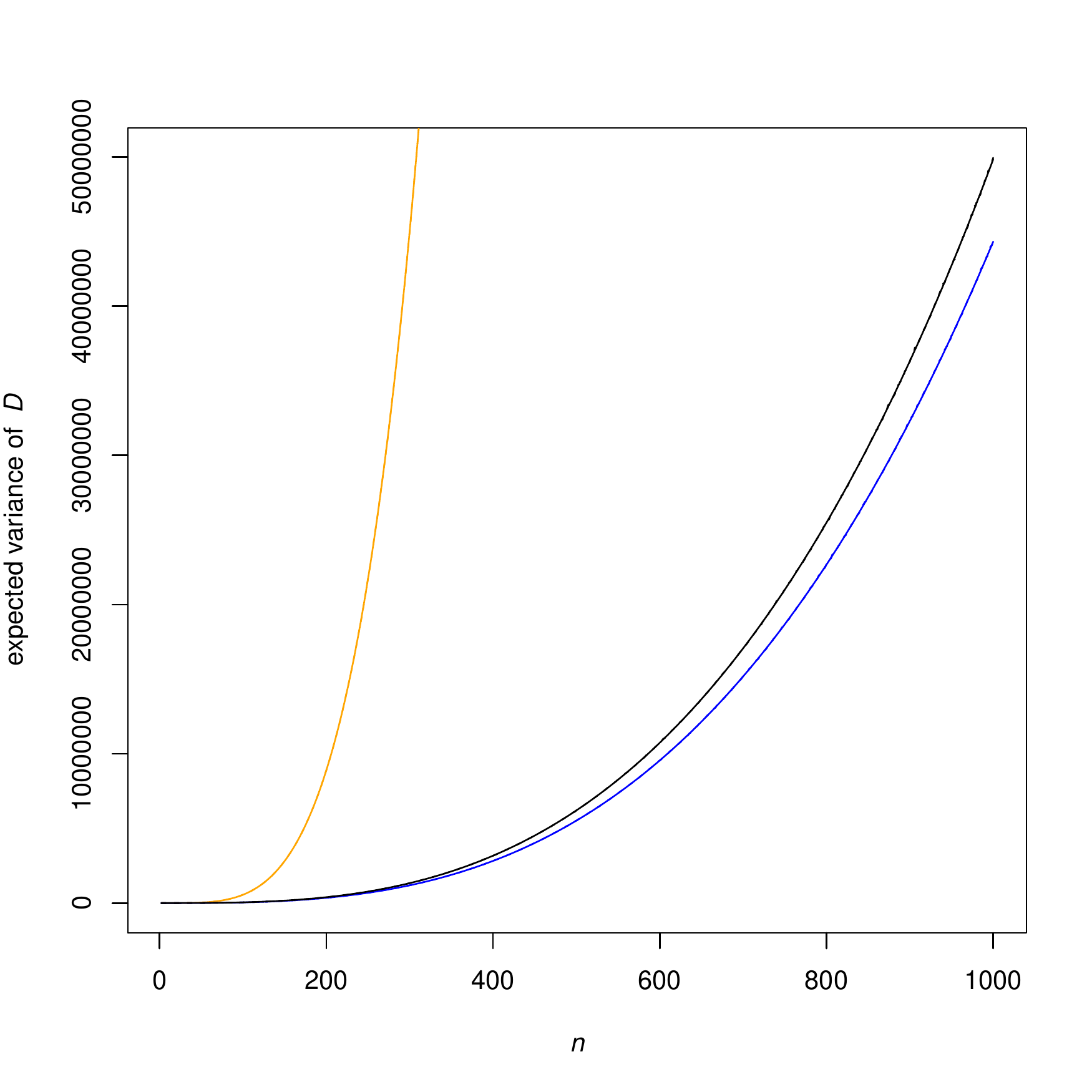}
\end{center}
\caption{\label{scaling_of_variance_random_trees_equation}  $\E_{rlt}\left[\V_{rla}[D]\right]$, the expected variance of $D$ in uniformly random labelled trees, as a function of $n$, the number of vertices of the tree (black). For reference,  $\V_{rla}[D]$ of a linear tree (blue) and a star tree (orange) are also shown. Solid lines indicate theoretical values while dashed lines indicate simulations. Dashed lines can hardly be seen due to the high accuracy of the computer estimations. }
\end{figure}
 
\section{Upper bounds}

\label{upper_bounds_section}

Here we aim to derive some elementary upper bounds of $D$ to help to show possibilities for future research in Section \ref{applications_section}. Additionally, this will complement the understanding of the variation of $D$ via $\V_{rla}[D]$ above. Suppose that $D_{max}$ is the maximum value of $D$ over the $n!$ linear arrangements of a graph. 
Trivially, 
\begin{equation*}
D_{max} \leq D(K_{n}).
\end{equation*}
However, this inequality is not very useful because $D(K_{n})$ depends only on $n$. Better upper bounds of $D$ for a given network can be obtained taking into account $m$. Obviously,
\begin{equation}
D_{max} \leq m(n - 1)^2,
\label{naive_sum_of_edge_lengths_max_equation}
\end{equation} 
as the maximum length of an edge is $n - 1$. We will derive upper bounds of $D_{max}$ applying some of the methods that Petit applied for deriving lower bounds of $D_{min}$ \cite{Petit2003a}. 

A tighter upper bound can be obtained with an analog of Petit's degree method \cite{Petit2003a} and the fact that $D$ can be defined equivalently as
\begin{equation}
D = \frac{1}{2} D_i,
\label{alternative_sum_of_edge_lengths_equation}
\end{equation}
where $D_i$ is the sum of the lengths of the edges involving the $i$-th vertex.
An upper bound of $D_i$ is given by placing $i$ at one end of the linear arrangement and its $k_i$ adjacent vertices as far as possible, i.e.   
\begin{equation}
D_i \leq \sum_{j=1}^{k_i} (n-j) = k_i \left(n - \frac{k_i + 1}{2}\right).
\label{sum_of_edge_length_of_vertex_max_equation}
\end{equation}
When $i$ is the hub of a star tree, $k_i = n - 1$ and then the right hand side of Eq. \ref{sum_of_edge_length_of_vertex_max_equation} becomes
\begin{equation}
{n \choose 2},
\end{equation}
matching Eq. \ref{sum_of_edge_length_star_tree_equation} with $\tau = 1$ of $\tau = n$. 

Eq. \ref{sum_of_edge_length_of_vertex_max_equation} transforms Eq. \ref{alternative_sum_of_edge_lengths_equation} into an upper bound of $D_{max}$ via the degree method ($DM$), i.e.  
\begin{equation}
D_{upper}^{DM} = \frac{1}{2}\left[ \left(n-\frac{1}{2}\right) \sum_{i=1}^n k_i - \frac{1}{2}\sum_{i=1}^n k_i^2 \right].
\end{equation}
Finally, the handshaking lemma \cite[p. 4]{Bollobas1998a} and the definition of $\left<k^2\right>$ give
\begin{equation}
D_{upper}^{DM} = m\left(n - \frac{1}{2}\right) - \frac{n}{4}\left<k^2\right>,
\label{sum_of_edge_lengths_max1_equation}
\end{equation}
which is an obvious improvement over Eq. \ref{naive_sum_of_edge_lengths_max_equation}. However, such an upper bound of $D$ performs poorly when applied to a complete graph, where $D$ is constant ($D_{max}=D_{min}$).  Applying Eq. \ref{sum_of_edge_lengths_max1_equation} to a complete graph of $n$ vertices, where 
\begin{eqnarray}
m(\completegraph) = {n \choose 2}\\
\left<k^2\right>(\completegraph) = (n-1)^2, 
\end{eqnarray}
routine calculations give 
\begin{equation*}
D_{max} \leq \frac{1}{4}(n-1)n^2,
\end{equation*} 
that is far from the real value of $D_{max}$ in Eq. \ref{sum_of_edge_lengths_complete_graph_equation}.

It is possible to get a general upper bound of $D_{max}$ that will match the true value of $D_{max}$ when applied to a complete graph. The method consist of noting that there are $n - d$ edges of length $d$ (for $1 \leq d \leq n - 1$) as in Petit's edges method \cite{Petit2003a}.
We define $F(d_0)$ as the maximum number of edges that can be formed of length within $[d_0, n-1]$, where $n-1$ is the maximum edge length. Then
\begin{eqnarray}
F(d_0) & = & \sum_{d = d_0}^{n - 1} n - d \nonumber \\
       & = & \sum_{i=1}^{n - d_0} i \nonumber \\
       & = & \frac{1}{2}(n - d_0)(n - d_0 + 1) \label{as_long_as_possible_edges_equation}
\end{eqnarray}
for $d_0 \in [0, n]$. 
We define $d_*$ as the smallest value of $d_0$ such that $F(d_0) \leq m$. 
Then we can build a network that maximizes $D$ for a given $m$ by taking all edges of length between $d_*$ and $n-1$, namely $F(d_*)$ edges, and $m - F(d_*)$ edges of length $d_* - 1$.
The sum of edge lengths of such a graph gives an upper bound of $D_{max}$, via the edges method ($EM$), i.e. 
\begin{eqnarray}
D_{upper}^{EM} & = & (m - F(d_*))(d_* - 1)  + \sum_{d=d_*}^{n-1} (n - d) d \nonumber \\
   & = & (m - F(d_*))(d_* - 1) \nonumber \\
   &   & + \frac{1}{6}(n- d_*) (n^2 + (n+3)d_* - 2d_*^2 - 1). \label{sum_of_edge_lengths_max_alternative_equation}
\end{eqnarray}
We want to derive $d_*$. 
Solving the equation $F(d_0) = m$, namely
\begin{equation*}
d_0^2 - (2n+1)d_0 + n(n+1) - 2m = 0
\end{equation*}
one finds two solutions for $d_0$, i.e.
\begin{equation*}
n + \frac{1}{2} \pm \frac{1}{2}\sqrt{8m + 1}
\end{equation*}   
and then
\begin{equation}
d_* = \left\lceil n + \frac{1}{2} - \frac{1}{2}\sqrt{8m + 1} \right\rceil. \label{critical_distance_equation}
\end{equation}
It is easy to check that $D_{upper}^{EM}$ (Eq. \ref{sum_of_edge_lengths_max_alternative_equation}) gives the actual maximum sum of edges lengths for extreme values of $m$. 
When $m = 0$, Eq. \ref{critical_distance_equation} gives $d_* = n$, $F(d_*)=0$ (Eq. \ref{as_long_as_possible_edges_equation}) and then Eq. \ref{sum_of_edge_lengths_max_alternative_equation} gives $D_{upper}^{EM} = 0$.
When $m = m(K_{n})$, Eq. \ref{critical_distance_equation} gives 
\begin{eqnarray*}
d_* & = & \left\lceil n + \frac{1}{2} - \frac{1}{2} \sqrt{8 {n \choose 2} + 1} \right\rceil \\
    & = & \left\lceil n + \frac{1}{2} - \frac{1}{2} \sqrt{(2n -1)^2} \right\rceil \\
    & = & 1
\end{eqnarray*} 
and also (recall Eq. \ref{as_long_as_possible_edges_equation})
\begin{equation}
F(d_*)= {n \choose 2}.  
\end{equation} 
Therefore, Eq. \ref{sum_of_edge_lengths_max_alternative_equation} gives $D_{upper}^{EM} = D(\completegraph)$ (Eq. \ref{sum_of_edge_lengths_complete_graph_equation}).

Taking the tightest of the upper bounds of $D_{max}$ that have been derived above (Eq. \ref{sum_of_edge_lengths_max_alternative_equation} and \ref{sum_of_edge_lengths_max1_equation}), one obtains a general upper bound of $D_{max}$, i.e. 
\begin{equation}
D_{upper} = \min (D_{upper}^{DM}, D_{upper}^{EM}).
\label{sum_of_edge_lengths_max_equation}
\end{equation}

Notice that $D_{upper}^{EM}$ can be calculated in constant time. In his pioneering research, Petit applied the edges method to derive a linear time algorithm to calculate a lower bound for $D_{min}$ \cite{Petit2003a}. The procedure of our derivation of $D_{upper}^{EM}$ could be applied to calculate a lower bound of $D_{min}$ in constant time. In particular, 
\begin{equation}
D_{min} \geq D(\completegraph) - D_{upper}^{EM},
\end{equation} 
where $D_{upper}^{EM}$ is calculated on a graph with $m(\completegraph) - m$ edges. 

\section{Applications}

\label{applications_section}

Our theoretical results on the distribution of $D$ in random linear arrangements have many potential applications. Here we only sketch some hoping that they illustrate the importance of our theoretical work and inspire future research.

\subsection{$z$-scoring of edge lengths}

\label{z_scoring_subsection}

The target of this application is the calculation of the mean edge length over a collection of networks, e.g., a treebank. A treebank is a collection of syntactic dependency trees as that of Fig. \ref{syntactic_dependency_tree_figure} \cite{Liu2008a}. 
The mean edge length of a network is defined as
\begin{eqnarray*}
\left<d \right> & = & D/m \\
                & = & \frac{1}{m} \sum_{i=1}^m d_i.
\end{eqnarray*}
Suppose a collection of $T$ networks where the $i$-th network has $n_i$ vertices, $m_i$ edges and $d_{ij}$ is the length of the $j$-th edge of the $i$-th network. The mean edge length of the collection can be defined as
\begin{eqnarray}
\left<d \right> = \frac{1}{M} \sum_{i=1}^T \sum_{j=1}^{m_i} d_{ij},
\label{mean_edge_length_in_collection_equation}
\end{eqnarray}
where
\begin{equation*}
M = \sum_{i=1}^T m_i
\end{equation*} 
is the total number of edges. 
If the networks are trees, then $m_i = n_i - 1$ and the mean edge length becomes
\begin{equation}
\left<d \right> = \frac{1}{N - T} \sum_{i=1}^T \sum_{j=1}^{m_i} d_{ij},
\label{mean_edge_length_in_collection_trees_equation}
\end{equation}
where 
\begin{equation*}
N = \sum_{i=1}^T n_i 
\end{equation*}
is the total number of vertices.
Eq. \ref{mean_edge_length_in_collection_trees_equation} matches the average edge length defined by Liu on collections of syntactic dependency trees. 
A general problem of Eq. \ref{mean_edge_length_in_collection_equation} is that the distribution of the inner summation, i.e.
\begin{equation*}
\sum_{j=1}^{m_i} d_{ij}, 
\end{equation*} 
depends on $n$, $m$ and $\left<k^2\right>$ under the null hypothesis (recall Eq. \ref{variance_of_sum_of_edge_lengths_equation}). Put differently, the mean edge length of the collection mixes lengths that may have different distributions under the null hypothesis. A $z$-score is a way to normalize the individual lengths to turn them more comparable. A $z$-score is a transformation of a random variable so that it has zero mean and unit standard deviation with respect to a certain distribution \cite{Kreyszig1979a}. Thanks to the theoretical results of our article, we can define a mean edge length over $z$-scores. First, notice that Eq. \ref{mean_edge_length_in_collection_equation} can be defined equivalently as
\begin{equation*}
\left<d \right> = \frac{1}{M} \sum_{i=1}^T D_i,
\end{equation*}
where $D_i$ is the sum of edge lengths of the $i$-th network. i.e.
\begin{equation}
D_i = \sum_{j=1}^{m_i} d_{ij}.
\end{equation} 
The mean $z$-scored edge length of the collection is
\begin{equation*}
\left<d \right>_z = \frac{1}{M} \sum_{i=1}^T z_i,
\end{equation*}
where $z_i$ is the $z$-score of the sum of edge lengths of the $i$-th network, i.e.  
\begin{equation}
z_i = \frac{D_i - \E_{rla}[D]_i}{\V_{rla}[D]_i^{1/2}},
\label{z_score_equation}
\end{equation}
where $\E_{rla}[D]_i$ and $\V_{rla}[D]_i$ are calculated applying the values of $n$, $m$ and $\left<k^2\right>$ of the $i$-th network to Eqs. \ref{1st_moment_of_sum_of_edge_lengths_equation} and \ref{variance_of_sum_of_edge_lengths_equation}.

Table \ref{syntactic_dependency_tree_table} allows one to calculate easily a $z$-scored value of $D$ for the network in Fig. \ref{syntactic_dependency_tree_figure}, 
i.e.
\begin{equation}
z = \frac{-56}{\sqrt{216.8}} \approx -3.803.
\label{z_score_example_equation}
\end{equation}
We hope that our outline stimulates further theoretical and empirical research on the problem of dependency distance normalization \cite{Lei2018a}. 

\begin{table}
\caption{\label{syntactic_dependency_tree_table} Summary of the statistical features of the network in Fig. \ref{syntactic_dependency_tree_figure}. }
\begin{indented}
\item[]
\begin{tabular}{@{}lr}
\br
Feature & Value \\
\mr
$n$ & 17 \\
$m$ & 16 \\
$D$ & 40 \\
$\E_{rla}[D]$ & 96 \\
$\left<k^2 \right>$ & $\frac{88}{17}$ \\
$\V_{rla}[D]$ & $\frac{1084}{5} = 216.8$ \\
\br
\end{tabular}
\end{indented}
\end{table}

\subsection{A test of significance of $D$}

\label{test_of_significance_subsection}

The aim of this application is a simple and fast test of whether $D$ is significantly small. In syntactic dependency trees, it has been found that $D$ is below $\E_{rla}[D]$ in general and that fact has been attributed to a general principle of edge length minimization \cite{Liu2017a}. To test that the value of $D$ of a real network is significantly low one uses $D_{rla}$, the value of $D$ of a random linear arrangement of the same network for reference. In particular, one has to show that the $P(D_{rla} \leq D)$, the probability that a random linear arrangement gives the same or a smaller value of $D$, is smaller than a certain significance level $\alpha$. One could calculate $P(D_{rla} \leq D)$ by brute force, as the proportion of the $n!$ permutations of the order of the vertices where $D_{rla} \leq D$. As this procedure is computationally unaffordable for sufficiently large $n$, it is convenient to use a Monte Carlo procedure to avoid the time consuming task of generating the $n!$ possible orderings of the vertices. In that procedure, one generates only $R$ uniformly random permutations and estimates $P(D_{rla} \leq D)$ as the proportion of the $R$ uniformly random permutations where $D_{rla} \leq D$. However, that Monte Carlo test is still time consuming if the network is large ($n$ is large) or a large $R$ is needed for accuracy. 
An alternative is to use well known inequalities that yield an upper bound of $P(D_{rla} \leq D)$ with little computational effort. To show the potential of this method, we chose a one-sided Chebychev inequality, also known as Cantelli's inequality, that for a random variable $x$ with expectation $\mu$ and standard deviation $\sigma$ gives \cite{Padulo2011a,Cantelli1910a} 
\begin{equation}
P(x - \mu \leq c \sigma) \leq \frac{1}{1+ c^2},
\end{equation}
where $c$ is a positive real number. 
Replacing $x$ by $D_{max} - D$ one obtains 
\begin{eqnarray*}
\mu = \E_{rla}[D_{max} - D] = D_{max} - \E_{rla}[D] \\
\sigma = \V_{rla}[D_{max} - D]^{1/2} = \V_{rla}[D]^{1/2}
\end{eqnarray*}
and finally
\begin{equation}
P\left(\E_{rla}[D] - D \leq c \V_{rla}[D]^{1/2}\right) \leq \frac{1}{1 + c^2}
\label{one_sided_Chebychev_equation}
\end{equation}
with a critical value of $c^*$ that is  
\begin{equation*}
c^* = \frac{\E_{rla}[D] - D}{\V_{rla}[D]^{1/2}}.
\end{equation*}
Notice that $-c^*$ is a $z$-score (recall Eq. \ref{z_score_equation}) and then Eq. \ref{z_score_example_equation}
gives
\begin{equation}
c^* =\frac{56}{216.8^{1/2}} \approx 3.803 \label{standard_deviations_equation} 
\end{equation} 
for the network in Fig. \ref{syntactic_dependency_tree_figure}.
Applying Eq. \ref{standard_deviations_equation} to  Eq. \ref{one_sided_Chebychev_equation} one obtains
\begin{equation*}
P(D_{rla} \leq D) \leq \frac{271}{4191} \approx 0.065.
\end{equation*}    
Thus, the one-sided Chebychev inequality does not support edge length minimization at a significance level of $\alpha = 0.05$. However, assuming the distribution of $x$ is symmetrical and unimodal, it follows that \cite{Padulo2011a, Popescu2005a}
\begin{equation*}
P(x - \mu \leq c \sigma) 
   \leq \left\{
      \begin{array}{ll}
         \frac{2}{9c^2} & \mbox{if~} c \geq 2/3 \\
         \frac{1}{2} & \mbox{if~} c \leq 2/3.
      \end{array}  
   \right.  
\end{equation*}
The substitution $x = D_{max} - D$ gives 
\begin{equation*}
P(\E_{rla}[D] - D \leq c \V_{rla}[D]^{1/2}) 
   \leq \left\{
      \begin{array}{ll}
         \frac{2}{9c^2} & \mbox{if~} c \geq 2/3 \\
         \frac{1}{2} & \mbox{if~} c \leq 2/3,
      \end{array}  
   \right.  
\end{equation*}
that is the version of the inequality needed in our application. 
Applying previous results for $c$ (Eq. \ref{standard_deviations_equation}), one obtains 
\begin{eqnarray*}
P(D_{rla} \leq D) & \leq & \frac{2}{9c^2} \\
                  & =    & \frac{271}{17640} \approx 0.015.
\end{eqnarray*}
If the distribution of $D_{rla}$ was actually symmetric and unimodal, we would conclude that $D$ is significantly small having spent a 	small amount of computational resources. This example illustrates the importance of further research on the properties of the distribution of $D_{rla}$. Notice that symmetricity and unimodality are just examples of properties for further research. It would also be possible to bound 
$P(D_{rla} \leq D)$ with the help of higher moments about zero, in particular, $\E_{rla}[D^3]$ and $\E_{rla}[D^4]$ \cite{Bhattacharyya1987a}.  

\subsection{Minimum linear arrangement problem}

\label{mla_subsection}

Here we wish to outline a potential contribution to the computationally hard problem of calculating $D_{min}$, known as the minimum linear arrangement problem in computer science \cite{Petit2003a}. In particular, our results may allow one to derive upper bounds of 
$D_{min}$ that can be used as random baselines to evaluate computational methods to calculate $D_{min}$ approximately \cite{Petit2003a,Caprara2011a}. 
A straightforward upper bound is obtained from general properties of expectation, namely, \cite[p. 188]{DeGroot1989a}
\begin{equation*}
D_{min} \leq \E_{rla}[D]
\end{equation*}
with equality if and only if probability mass is concentrated on $D_{min}$ ($P(D = D_{min})=1$, $P(D > D_{min})=0$). Such a baseline is known as the random layout \cite{Petit2003a}.
Here we will derive an upper bound of $D_{min}$ of the form
\begin{equation}
D_{min} \leq \E_{rla}[D] - \Delta,
\label{upper_bound_of_D_min_equation}
\end{equation}
where $\Delta$ is a positive quantity that is a function of properties of the distribution of $D_{rla}$. Our target are bounds that are of low computational cost.

Bathia-Davis' inequality bounds variance above based on $\E_{rla}[D]$, $D_{min}$ and $D_{max}$, namely the average, the minimum and the maximum value of $D$ over the $n!$ linear arrangements of a given network. In particular, this inequality states that \cite{Bathia2000a} 
\begin{equation*}
\V_{rla}[D] \leq (D_{max} - \E_{rla}[D])(\E_{rla}[D] - D_{min})
\end{equation*}  
and is equivalent to Eq. \ref{upper_bound_of_D_min_equation} with
\begin{equation*}
\Delta = \frac{\V_{rla}[D]}{D_{max} - \E_{rla}[D]}
\end{equation*}
assuming $D_{max} \neq \E_{rla}[D]$.
$\E_{rla}[D]$ and $\V_{rla}[D]$ can be calculated easily via Eqs. \ref{1st_moment_of_sum_of_edge_lengths_equation} and \ref{variance_of_sum_of_edge_lengths_equation}. $D_{max}$ can be obtained with an algorithm to solve the maximum linear arrangement problem but the calculation is computationally expensive for arbitrary networks \cite{Hassin2001a}. Sacrificing accuracy of the upper bound of $D_{min}$ for the sake of computational efficiency, $D_{max}$ can be replaced by $D_{upper}$, the upper bound of $D_{max}$ obtained in Section \ref{upper_bounds_section} (Eq. \ref{sum_of_edge_lengths_max_equation}). 

Sherma et al's inequality introduces $\W_{rla}[D]$, the third central moment of $D$ in a random linear arrangement, namely 
\begin{equation}
\W_{rla}[D] = \E_{rla}[D^3] + 2\E_{rla}[D]^3 -  3\E_{rla}[D]\E_{rla}[D^2].
\label{third_central_moment_equation}
\end{equation}
In particular, Sherma et al's inequality states that \cite{Sharma2010a}
\begin{equation*}
\V_{rla}[D] \leq \frac{1}{4} (D_{max} - D_{min})^2 - \left(\frac{\W_{rla}[D]}{2\V_{rla}[D]} \right)^2,
\end{equation*}
giving
\begin{equation*}
D_{min } \leq D_{max} - 2 \left(\V_{rla}[D] + \frac{\W_{rla}[D]}{2\V_{rla}[D]} \right)^{1/2}.
\end{equation*}
The last result and the definition of $\W_{rla}$ above (Eq. \ref{third_central_moment_equation}) 
show the importance of investigating further the properties of the distribution of $D$ (e.g., $\E_{rla}[D^3]$ and $D_{max}$). 

\section{Discussion}

\label{discussion_section}

Throughout this article, we have deepened our understanding of the distribution of $D$. On the one hand, 
we have presented two different derivations of the $\E_\phi$'s (Section \ref{second_moment_section} and \ref{alternative_derivation_section}), obtaining compact formulae for $\E_{rla}[D^2]$ and $\V_{rla}[D]$ in general (Section \ref{second_moment_section}) and for specific networks (Section \ref{complete_graph_subsection} and \ref{validation_section}). Table \ref{summary_table} summarizes all the values of $\E_{rla}[D^2]$ and $\V_{rla}[D]$ that have been obtained for specific networks with given $n$ throughout this article or recycled from previous work. 
On the other hand, we have have obtained upper bounds of $D_{max}$ (Section \ref{upper_bounds_section}) and suggested new tracks for exploring upper bounds on $D_{min}$ that may stimulate further research in computer science \cite{Petit2011a} (Section \ref{mla_subsection}).

\begin{table}
\caption{\label{summary_table} Summary of 1st and 2nd moment about zero of $D$, the sum of edge lengths, and the variance of $D$, in random linear arrangements for relevant networks in this article. By default, equations are valid for $n \geq 0$. Equations valid for $n\geq 1$ are marked with $^*$. Those valid for $n\geq 2$ are marked with $^{**}$.} 
\begin{indented}
\item[]
\begin{tabular}{@{}llrr}
\br
Graph   & $\E_{rla}[D]$ & $\E_{rla}[D^2]$ & $\V_{rla}[D]$ \\ % \hline
\mr
$m = 0$              & 0 & 0 & 0 \\
$m = 1$              & $\frac{1}{3}(n+1)^{**}$ & $\frac{1}{6}n(n+1)$ & $\frac{1}{18}(n+1)(n-2)$ \\
$\lineartree$        & $\frac{1}{3}(n^2 - 1)^*$ & $\frac{1}{90} (n + 1) (10n^3 - 6n^2 - 25n +24)^{**}$ & $\frac{1}{90} (n + 1) (n-2) (4n - 7)^{**}$ \\
$\startree$          & $\frac{1}{3}(n^2 - 1)^*$ & $\frac{1}{60}(n^2 - 1)(7n^2-8)^*$ & $\frac{1}{180}(n^2 -1)(n^2 - 4)^*$ \\
$m = {n \choose 2}$  & $\frac{1}{6}(n^2 - 1)n$ & $\left[ \frac{1}{6}(n^2-1)n \right]^2$ & 0 \\
\br
\end{tabular}
\end{indented}
\end{table}

We have applied the theoretical results on the variance of $D$ to a couple of network ensembles. In 
Erd\H{o}s-R\'enyi graphs, we have found that the expected $\V_{rla}[D]$ as a function of $m$ evolves following a bell-shape peaking when the density of links is about 1/2 while $n$ remains constant (Section \ref{Erdos_Renyi_graph_section}). In uniformly random labelled trees, we have found that the expected $\V_{rla}[D]$ as a function of $n$ scales asymptotically following a power-law of $n$ (Section \ref{uniformly_random_labelled_trees_section}). Other classes of random networks with more realistic characteristics should be investigated \cite{Newman2010a,Ferrer2017a}.
In addition, we have applied the theoretical results to obtain $z$-scored measures of edge length (Section \ref{z_scoring_subsection}) and to develop a simple test of significance of $D$ that can be very helpful in language research \cite{Liu2017a} (Section \ref{test_of_significance_subsection}).
We hope that our work stimulates further research on the properties of the distribution of $D$ in random linear arrangements and further applications across disciplines. 
 
\ack
We are grateful to L. Alemany-Puig, C. Card\'o, N. Catal\`a and an anonymous reviewer for valuable comments and helping to improve the quality of this manuscript. We also thank Jordi Petit for helpful discussions on the minimum linear arrangement problem. This research was supported by the grant TIN2017-89244-R from MINECO (Ministerio de Economia, Industria y Competitividad) and the recognition 2017SGR-856
(MACDA) from AGAUR (Generalitat de Catalunya).

\appendix

\section{Validation}

\label{validation_section}

The compact formula for $\E_0$ in Eq. \ref{expectation_type_0_equation} is validated using a computational procedure that checks that $\E_0$ matches the average of $d_i d_j$ for a pair of edges of type 0 over the $n!$ linear arrangements of the two edges for a given $n$.   
The same validation procedure is applied to the compact formulae for $\E_1$ in Eq. \ref{expectation_type_1_equation} and the compact formulae for $\E_2$ in Eq. \ref{expectation_type_2_equation}. On top of it, each $E_\phi$ is tested for $n \in [4-\phi, 15]$.  

The formulae for $\V_{rla}[D]$ and $\E_{rla}[D^2]$ are validated with the help of graphs where $\V_{rla}[D]$ is known {\em a priori} or easy to derive independently:
\begin{itemize}
\item
A graph with minimum $m$, namely $m=0$, and then $\V_{rla}[D]=\E_{rla}[D^2]=\E_{rla}[D] = \left<k^2 \right>= 0$. Checking this is straightforward with the help of Eq. \ref{1st_moment_of_sum_of_edge_lengths_equation} and \ref{2nd_moment_of_sum_of_edge_lengths_equation}. 
\item
A graph with maximum $m$, namely a complete graph, where $\V_{rla}[D]=0$ because $D$ is constant. This implies that $\E_{rla}[D^2] = \E_{rla}[D]^2$. Applying the definition of $m(\completegraph)$ (Eq. \ref{edges_complete_graph_equation}) and 
$\left<k^2\right>(\completegraph) = (n-1)^2$ to Eq. \ref{2nd_moment_of_sum_of_edge_lengths_equation}, one obtains the definition of $\E_{rla}[D(\completegraph)^2]$ in Eq. \ref{2nd_moment_of_sum_of_edge_lengths_complete_graph_equation} as expected. 
\item
A graph with $m=1$. In this case, $D = d$, where $d$ is the length of the single edge in a uniformly random linear arrangement and then $\V_{rla}[D] = \V_{rla}[d]$ and $\E_{rla}[D^2] = \E_{rla}[d^2]$, that are provided in Eqs. \ref{2nd_moment_of_sum_of_edge_lengths_one_edge_equation} and 
\ref{variance_of_sum_of_edge_lengths_one_edge_equation}.
Eq. \ref{2nd_moment_of_sum_of_edge_lengths_equation} with $m=1$ and $\left<k^2\right>=2/n$ produces Eq. \ref{2nd_moment_of_sum_of_edge_lengths_one_edge_equation} as expected.
\item
In a star tree of $n$ vertices, i.e. $\startree$, $\V_{rla}[D]$ and $\E_{rla}[D^2]$ are easy to derive independently.
In such a tree, $D$ is determined by $\tau$, the position of the hub vertex in the linear arrangement ($1 \leq \tau  \leq n$). It is known that $D_\tau$, the value of $D$ as a function of $\tau$, is  \cite{Ferrer2013e} 
\begin{equation}
D_\tau = \tau^2 - (n+1)\tau + \frac{1}{2}n(n+1).
\label{sum_of_edge_length_star_tree_equation}
\end{equation}
It is easy to see that 
\begin{equation*}
\E_{rla}[D^2](\startree) = \frac{1}{n} \sum_{\tau=1}^n D_\tau^2
\end{equation*}
and then
\begin{eqnarray}
\E_{rla}[D^2] (\startree) = \frac{1}{60}(n+1)(n-1)(7n^2-8)
\label{2nd_moment_of_sum_of_edge_lengths_star_tree_equation}
\end{eqnarray}
after some algebra. % (comprovar aquesta formula via simulacio ???).
The combination of Eqs. \ref{1st_moment_of_sum_of_edge_lengths_equation} with $m = n-1$ and Eq. \ref{2nd_moment_of_sum_of_edge_lengths_star_tree_equation} gives
\begin{eqnarray}
\V_{rla}[D](\startree) & = & \frac{1}{60}(n+1)(n-1)(7n^2-8) - \left[\frac{1}{3}(n+1)(n-1)\right]^2 \nonumber \\
     & = & \frac{1}{180}(n+1)(n-1)(n+2)(n-2). \label{variance_of_sum_of_edge_lengths_star_tree_equation}
\end{eqnarray}
It is easy to check that Eq. \ref{2nd_moment_of_sum_of_edge_lengths_tree_equation} with $\left<k^2\right>(\startree)=n-1$ \cite{Ferrer2013b} gives \ref{2nd_moment_of_sum_of_edge_lengths_star_tree_equation} as expected. 
\end{itemize}

The test cases above cover only a small set of values of $m$, namely 
\begin{equation*}
m \in \left\{0,1, n-1, {n \choose 2} \right\}, 
\end{equation*}
missing many graphs for $m=n-1$. We wish to test the equations with graphs that cover all the possible values of $m$ given $n$, i.e. any $m$ such that  
\begin{equation*} 
m \in \left[0, {n \choose 2}\right]. 
\end{equation*}
This is satisfied in Section  \ref{Erdos_Renyi_graph_section} with the help of Erd\H{o}s-R\'enyi graphs (Fig. \ref{Erdos_Renyi_graph_figure}). Additional testing is performed on random trees in Section  \ref{uniformly_random_labelled_trees_section}.

\section{Alternative derivation}

\label{alternative_derivation_section}

We have derived formulae for $\E_{rla}[D^2]$ and $\V_{rla}[D]$ applying a concrete method: $\E_2$ is borrowed from previous work, $\E_1$ is derived independently and $\E_0$ is obtained as the solution of a linear equation of $\E_{rla}[D^2]$ on complete graphs, namely
\begin{equation*}
\E_{rla}[D^2](\completegraph) = f_0(\completegraph) \E_0 + f_1(\completegraph) \E_1 + f_2(\completegraph) \E_2,
\end{equation*}
where $\E_0$ is the only unknown ($\E_{rla}[D(\completegraph)^2]$ and $f_\phi(\completegraph)$ for $\phi \in [0,2]$ have been obtained independently). 

Alternatively, we could have derived $\E_0$ and $\E_1$ from scratch as the solutions of a system of two linear equations, i.e. 
\begin{eqnarray}
\E_{rla}[D^2](\completegraph) = f_0(\completegraph) \E_0 + f_1(\completegraph) \E_1 + f_2(\completegraph) \E_2 \label{1st_linear_equation}\\ 
\E_{rla}[D^2](\startree) = f_0(\startree) \E_0 + f_1(\startree) \E_1 + f_2(\startree) \E_2,  \label{2nd_linear_equation} 
\end{eqnarray}
where $\E_0$ and $\E_1$ are the only unknowns. In a star tree, 
\begin{itemize}
\item 
$f_2(\startree) = m = n - 1$ for $n\geq 1$.
\item 
$f_0(\startree) = 2q = 0$ because all pairs of edges share the hub vertex.
\item 
Combining  
\begin{equation*}
(n - 1)^2 = f_0(\startree) + f_2(\startree) + f_1 (\startree)
\end{equation*} 
and the values of $f_0(\startree)$ and $f_2(\startree)$ above, one obtains
\begin{equation*}
f_1 (\startree) = (n - 1)(n - 2).
\end{equation*}
\end{itemize}
Now we have all the information that is needed to solve the system of two linear equations. 
The second linear equation (Eq. \ref{2nd_linear_equation}) gives
\begin{equation*}
\E_1 = \frac{\E_{rla}[D^2](\startree)  -  (n - 1) \E_2}{(n-1)(n-2)}.    
\end{equation*}
Applying Eqs. \ref{2nd_moment_of_sum_of_edge_lengths_star_tree_equation} and \ref{2nd_moment_of_sum_of_edge_lengths_one_edge_equation}, one recovers the definition of $\E_1$ in Eq. \ref{expectation_type_1_equation}.
Finally notice that the first linear equation (Eq. \ref{1st_linear_equation}) gives
\begin{equation*}
\E_0 = \frac{\E_{rla}[D^2](\completegraph)  - f_1(\completegraph) \E_1 - f_2(\completegraph) \E_2}{f_0(\completegraph)},
\end{equation*}
that corresponds to Eq. \ref{precursor_equation}, that we have already used to obtain $\E_0$ in Section \ref{second_moment_section}. 
% Finally, notice that we have already shown that Eq. \ref{expectation_type_0_equation} follows from applying Eq. \ref{expectation_type_1_equation} to the 1st linear equation (Eq. \ref{1st_linear_equation}).

\section*{References}

\bibliographystyle{unsrt}

% \bibliography{../../../../Dropbox/biblio/rferrericancho,../../../../Dropbox/biblio/complex,../../../../Dropbox/biblio/ling,../../../../Dropbox/biblio/cl,../../../../Dropbox/biblio/cs,../../../../Dropbox/biblio/maths}

\begin{thebibliography}{10}

\bibitem{Barthelemy2011a}
M.~Barth\'elemy.
\newblock Spatial networks.
\newblock {\em Physics Reports}, 499(1):1 -- 101, 2011.

\bibitem{Krioukov2010a}
D.~Krioukov, F.~Papadopoulos, M.~Kitsak, A.~Vahdat, and M.~Bogu{\~n}\'a.
\newblock Hyperbolic geometry of complex networks.
\newblock {\em Physical Review E}, 82:036106, 2010.

\bibitem{Ferrer2017c}
R.~{Ferrer-i-Cancho}.
\newblock Towards a theory of word order. {Comment} on "{Dependency} distance:
  a new perspective on syntactic patterns in natural language" by {Haitao Liu}
  et al.
\newblock {\em Physics of Life Reviews}, 21:218--220, 2014.

\bibitem{Bullmore2009a}
E.~Bullmore and O.~Sporns.
\newblock Complex brain networks: graph theoretical analysis of structural and
  functional systems.
\newblock {\em Nature Reviews Neuroscience}, 10:186--198, 2009.

\bibitem{Ferrer2004b}
R.~{Ferrer-i-Cancho}.
\newblock {Euclidean} distance between syntactically linked words.
\newblock {\em Physical Review E}, 70:056135, 2004.

\bibitem{Barthelemy2018a}
M.~Barth\'elemy.
\newblock {\em Morphogenesis of Spatial Networks}.
\newblock Springer, Cham, 2018.

\bibitem{Ferrer2016a}
R.~{Ferrer-i-Cancho} and C.~G\'omez-Rodr\'iguez.
\newblock Liberating language research from dogmas of the 20th century.
\newblock {\em Glottometrics}, 33:33--34, 2016.

\bibitem{Liu2017a}
H.~Liu, C.~Xu, and J.~Liang.
\newblock Dependency distance: a new perspective on syntactic patterns in
  natural languages.
\newblock {\em Physics of Life Reviews}, 21:171--193, 2017.

\bibitem{Diaz2002}
J.~D\'iaz, J.~Petit, and M.~Serna.
\newblock A survey of graph layout problems.
\newblock {\em ACM Computing Surveys}, 34:313--356, 2002.

\bibitem{Petit2011a}
J.~Petit.
\newblock Addenda to the survey of layout problems.
\newblock {\em Bulletin of the European Association for Theoretical Computer
  Science}, 105:177--201, 2011.

\bibitem{Esteban2016a}
J.~L. Esteban, R.~{Ferrer-i-Cancho}, and C.~G\'omez-Rodr\'iguez.
\newblock The scaling of the minimum sum of edge lengths in uniformly random
  trees.
\newblock {\em Journal of Statistical Mechanics}, page 063401, 2016.

\bibitem{Gomez2016a}
C.~G\'omez-Rodr\'iguez and R.~{Ferrer-i-Cancho}.
\newblock Scarcity of crossing dependencies: a direct outcome of a specific
  constraint?
\newblock {\em Physical Review E}, 96:062304, 2017.

\bibitem{Ferrer2013b}
R.~{Ferrer-i-Cancho}.
\newblock Hubiness, length, crossings and their relationships in dependency
  trees.
\newblock {\em Glottometrics}, 25:1--21, 2013.

\bibitem{Zornig1984a}
P.~Z\"{o}rnig.
\newblock The distribution of the distance between like elements in a sequence
  {I}.
\newblock {\em Glottometrika}, 6:1--15, 1984.

\bibitem{Ferrer2016d}
R.~{Ferrer-i-Cancho}.
\newblock Non-crossing dependencies: least effort, not grammar.
\newblock In A.~Mehler, A.~L{\"u}cking, S.~Banisch, P.~Blanchard, and B.~Job,
  editors, {\em Towards a theoretical framework for analyzing complex
  linguistic networks}, pages 203--234. Springer, Berlin, 2016.

\bibitem{Piazza1991a}
B.L. Piazza, R.D. Ringeisen, and S.K. Stueckle.
\newblock Properties of nonminimum crossings for some classes of graphs.
\newblock In Yousef~Alavi et~al., editor, {\em Proceedings of the 6th
  International Conference on Graph theory, combinatorics, and applications},
  volume~2, pages 975--989, New York, 1991. Wiley.

\bibitem{Bollobas1998a}
B.~Bollob\'as.
\newblock {\em Modern graph theory}.
\newblock Springer-Verlag, 1998.

\bibitem{Erdos1959a}
P.~Erd\H{o}s and A.~R\'enyi.
\newblock On random graphs {I}.
\newblock {\em Publicationes Mathematicae}, 6:290--297, 1959.

\bibitem{Bollobas2002a}
B.~Bollob\'as and O.~Riordan.
\newblock Mathematical results on scale-free random graphs.
\newblock In S.~Bornholdt and H.~Schuster, editors, {\em Handbook of graphs and
  networks: from the genome to the Internet}, pages 1--34. Wiley-VCH, Berlin,
  2003.

\bibitem{Erdos1960a}
P.~Erd\H{o}s and A.~R\'enyi.
\newblock On the evolution of random graphs.
\newblock {\em Publications of the Mathematical Institute of the Hungarian
  Academy of Sciences}, 5:17--61, 1960.

\bibitem{Liu2007a}
H.~Liu.
\newblock Probability distribution of dependency distance.
\newblock {\em Glottometrics}, 15:1--12, 2007.

\bibitem{Liu2008a}
H.~Liu.
\newblock Dependency distance as a metric of language comprehension difficulty.
\newblock {\em Journal of Cognitive Science}, 9:159--191, 2008.

\bibitem{Futrell2015a}
R.~Futrell, K.~Mahowald, and E.~Gibson.
\newblock Large-scale evidence of dependency length minimization in 37
  languages.
\newblock {\em Proceedings of the National Academy of Sciences USA},
  112(33):10336--10341, 2015.

\bibitem{Ouyang2017a}
J.~Ouyang and J.~Jiang.
\newblock Can the probability distribution of dependency distance measure
  language proficiency of second language learners?
\newblock {\em Journal of Quantitative Linguistics}, 25:295--313, 2018.

\bibitem{Lei2018a}
L.~Lei and M.~L. Jockers.
\newblock Normalized dependence distance: Proposing a new measure.
\newblock {\em Journal of Quantitative Linguistics}, 2018.

\bibitem{Ferrer2006d}
R.~{Ferrer-i-Cancho}.
\newblock Why do syntactic links not cross?
\newblock {\em Europhysics Letters}, 76(6):1228--1235, 2006.

\bibitem{Ferrer2014c}
R.~{Ferrer-i-Cancho}.
\newblock A stronger null hypothesis for crossing dependencies.
\newblock {\em Europhysics Letters}, 108:58003, 2014.

\bibitem{Aldous1990a}
D.~Aldous.
\newblock The random walk construction of uniform spanning trees and uniform
  labelled trees.
\newblock {\em SIAM J. Disc. Math.}, 3:450--465, 1990.

\bibitem{Broder1989a}
A.~Broder.
\newblock Generating random spanning trees.
\newblock In {\em Symp. Foundations of Computer Sci., IEEE}, pages 442--447,
  New York, 1989.

\bibitem{Cayley1889a}
A.~Cayley.
\newblock A theorem on trees.
\newblock {\em Quart. J. Math}, 23:376--378, 1889.

\bibitem{Moon1970a}
J.~Moon.
\newblock Counting labelled trees.
\newblock In {\em Canadian Math. Cong.}, 1970.

\bibitem{Noy1998a}
M.~Noy.
\newblock Enumeration of noncrossing trees on a circle.
\newblock {\em Discrete Mathematics}, 180:301--313, 1998.

\bibitem{Pruefer1918a}
H.~Pr\"ufer.
\newblock {Neuer Beweis eines Satzes \"uber Permutationen}.
\newblock {\em Arch. Math. Phys}, 27:742--744, 1918.

\bibitem{Alonso1995a}
L.~Alonso and R.~Schott.
\newblock {\em Random generation of trees. Random generators in computer
  science}.
\newblock Springer, Dordrecht, 1995.

\bibitem{Petit2003a}
J.~Petit.
\newblock Experiments on the minimum linear arrangement problem.
\newblock {\em Journal of Experimental Algorithmics}, 8, 2003.

\bibitem{Kreyszig1979a}
E.~Kreyszig.
\newblock {\em Advanced Engineering Mathematics}.
\newblock Wiley, 4th edition, 1979.

\bibitem{Padulo2011a}
M.~Padulo and M.~D. Guenov.
\newblock Worst-case robust design optimization under distributional
  assumptions.
\newblock {\em International Journal for Numerical Methods in Engineering},
  88(8):797--816, 2011.

\bibitem{Cantelli1910a}
F.~P. Cantelli.
\newblock Intorno ad un teorema fondamentale della teoria del rischio.
\newblock {\em Bollettino dell'Associazione degli Attuari Italiani}, 24:1--23,
  1910.

\bibitem{Popescu2005a}
I.~Popescu.
\newblock A semidefinite programming approach to optimal-moment bounds for
  convex classes of distributions.
\newblock {\em Mathematics of Operations Research}, 30(3):632–657, 2005.

\bibitem{Bhattacharyya1987a}
B.B. Bhattacharyya.
\newblock One sided {Chebyshev} inequality when the first four moments are
  known.
\newblock {\em Communications in Statistics - Theory and Methods},
  16(9):2789--2791, 1987.

\bibitem{Caprara2011a}
A.~Caprara, M.~Oswald, G.~Reinelt, R.~Schwarz, and E.~Traversi.
\newblock Optimal linear arrangements using betweenness variables.
\newblock {\em Math. Program. Comput.}, 3(3):261--280, 2011.

\bibitem{DeGroot1989a}
M.~H. DeGroot.
\newblock {\em Probability and statistics}.
\newblock Addison-Wesley, Reading, MA, 1989.
\newblock 2nd edition.

\bibitem{Bathia2000a}
R.~Bhatia and C.~Davis.
\newblock A better bound on the variance.
\newblock {\em The American Mathematical Monthly}, 107(4):353--357, 2000.

\bibitem{Hassin2001a}
R.~Hassin and S.~Rubinstein.
\newblock Approximation algorithms for maximum linear arrangement.
\newblock {\em Information Processing Letters}, 80(4):171 -- 177, 2001.

\bibitem{Sharma2010a}
R.~Sharma, M.~Gupta, and G.~Kapoor.
\newblock Some better bounds on variance with applications.
\newblock {\em Journal of Mathematical Inequalities}, 4(3), 2010.

\bibitem{Newman2010a}
M.~E.~J. Newman.
\newblock {\em Networks. An introduction}.
\newblock Oxford University Press, Oxford, 2010.

\bibitem{Ferrer2017a}
C.~{G\'omez-Rodr{\'i}guez} R.~{Ferrer-i-Cancho} and J.~L. Esteban.
\newblock Are crossing dependencies really scarce?
\newblock {\em Physica A}, 493:311--329, 2018.

\bibitem{Ferrer2013e}
R.~{Ferrer-i-Cancho}.
\newblock The placement of the head that minimizes online memory. {A} complex
  systems approach.
\newblock {\em Language Dynamics and Change}, 5:114--137, 2015.

\end{thebibliography}

\end{document}